\documentclass[review,3p]{IEEEtran}

\usepackage{mathrsfs,latexsym,float,enumerate,amssymb,amsmath,bm,color,lineno}
\usepackage{amscd,verbatim}
\usepackage{graphicx}
\usepackage[colorlinks=true,linkcolor=blue,citecolor=blue]{hyperref}
\usepackage{booktabs}
\usepackage{subfigure}

\newtheorem{theorem}{Theorem}[section]

\newtheorem{proposition}{Proposition}[section]

\newtheorem{definition}{Definition}[section]
\newtheorem{example}{Example}

\newtheorem{remark}{Remark}

\usepackage{pifont}

\begin{document}

\title{A Monotonous Intuitionistic Fuzzy TOPSIS Method under General Linear Orders
via Admissible Distance Measures}

\author{Xinxing Wu, Zhiyi Zhu, Chuan Chen, Guanrong Chen,~\IEEEmembership{Life Fellow,~IEEE,} Peide Liu 
\thanks{Manuscript received xx 00, 201x; revised xx 00, 201x.}
\thanks{X. Wu is with (1) the School of Sciences, Southwest Petroleum University, Chengdu, Sichuan 610500, China;
(2) the Institute for Artificial Intelligence, Southwest Petroleum University, Chengdu, Sichuan 610500, China.
(3) the Zhuhai College of Jilin University, Zhuhai, Guangdong 519041, China.
e-mail: (wuxinxing5201314@163.com).}
\thanks{Z. Zhu is with the School of Sciences, Southwest Petroleum University, Chengdu, Sichuan 610500, China.
e-mail: (zhuzhiyi2019@163.com).}
\thanks{C. Chen is with the School of Sciences, Southwest Petroleum University, Chengdu, Sichuan 610500, China.
e-mail: (chenchuan975201314@163.com).}
\thanks{G. Chen is with the Department of Electrical Engineering, City University of
Hong Kong, Hong Kong SAR, China.
e-mail: (eegchen@cityu.edu.hk).}
\thanks{P. Liu is with a School of Management Science and Engineering, Shandong University
of Finance and Economics, Jinan, Shandong 250014, China. e-mail: (peide.liu@gmail.com).}
\thanks{All correspondences should be addressed to X. Wu and P. Liu.}
\thanks{This work was supported by the National Natural Science Foundation of China
(Nos.~11601449 and 71771140), the Key Natural Science Foundation of Universities in
Guangdong Province (No. 2019KZDXM027), and the Natural Science Foundation of
Sichuan Province (No.~2022NSFSC1821).}
}


\maketitle

\begin{abstract}
All intuitionistic fuzzy TOPSIS methods contain two key elements: (1) the order structure, which can
affect the choices of positive ideal-points and negative ideal-points, and construction of admissible
distance/similarity measures; (2) the distance/similarity measure, which is closely related to the
values of the relative closeness degrees and determines the accuracy and rationality of decision-making.
For the order structure, many efforts are devoted to constructing some score functions, which can strictly
distinguish different intuitionistic fuzzy values (IFVs) and preserve the natural partial order for IFVs.
This paper proves that such a score function does not exist, namely the application of a single monotonous
and continuous function does not distinguish all IFVs. For the distance or similarity measure, some examples
are given to show that classical similarity measures based on the normalized Euclidean distance and normalized
Minkowski distance do not meet the axiomatic definition of intuitionistic fuzzy similarity measures. Moreover,
some illustrative examples are given to show that classical intuitionistic fuzzy TOPSIS methods do not ensure
the monotonicity with the natural partial order or linear orders, which may yield some counter-intuitive
results. To overcome the limitation of non-monotonicity, we propose a novel intuitionistic fuzzy TOPSIS method,
using three new admissible distances with the linear orders measured by a score degree/similarity
function and accuracy degree, or two aggregation functions, and prove that the proposed TOPSIS method is monotonous
under these three linear orders. This is the first result with a strict mathematical proof on the monotonicity
 with the linear orders for the intuitionistic fuzzy TOPSIS method. Finally, we show two practical examples and
comparative analysis with other decision-making methods to illustrate the efficiency of the developed TOPSIS method.
\end{abstract}

\begin{IEEEkeywords}
Intuitionistic fuzzy set, Distance measure, Similarity measure, TOPSIS,
Multi-attribute decision making.
\end{IEEEkeywords}

\IEEEpeerreviewmaketitle

\section{Introduction}
\IEEEPARstart{Z}{adeh} (1965)~\cite{Za1965} established the fuzzy set (FS) theory by applying
membership degree to measure the importance of a fuzzy element, which generalized
the concept of crisp set, characterized by a characteristic function taking value $0$ or
$1$, by taking any value in the unit interval $[0, 1]$. However, due to the limitation
of a membership function that only indicates two (supporting and opposing) opposite states
of fuzziness, the fuzzy set cannot express the neutral state of ``this and also that".
According to this, Atanassov (1986)~\cite{Ata1986} generalized Zadeh's fuzzy set
by proposing the concept of intuitionistic fuzzy sets (IFSs) (see also \cite{Ata1999}), characterized
by a membership function and a non-membership function meeting the condition
that the sum of the membership degree and the non-membership degree at every point
is less than or equal to $1$. Every pair of membership and non-membership degrees
for IFSs was called an intuitionistic fuzzy value (IFV) by Xu~\cite{Xu2007}.
Thereafter, various multi-attribute decision making (MADM) methods {under the
intuitionistic fuzzy framework} were developed and widely applied.
This paper focuses on the intuitionistic fuzzy TOPSIS method.

Being one of the most well-known MADM methods, TOPSIS was first proposed by Hwang and Yoon
(1981)~\cite{HY1981}. The main idea of the TOPSIS is that the most desirable
alternative should be nearest from the positive ideal-point and meanwhile furthest
from the negative ideal-point. Noting that the total order structure `$\leq$' and the
absolute distance `$|\cdot |$' on the real line $\mathbb{R}$ are admissible (the bigger
the real number, the nearer from the maximum, and the further from the minimum). This can
naturally guarantee the establishment of the TOPSIS and further ensure its monotonicity.
{Due to} the complex two-dimensional structure of the space of all IFVs, all existing IF distance
measurements for IFVs are not admissible under linear orders on IFVs. Therefore,
we~\cite{WWLCZ} introduced the concept of admissible distance measure with the
linear order `$\leq_{_{\mathrm{XY}}}$' of Xu and Yager~\cite{XY2006} and constructed
an admissible distance measure $\varrho$ with the linear order
`$\leq_{_{\mathrm{XY}}}$'.

Through careful analysis of the TOPSIS method, it is not difficult to find that
the TOPSIS method contains two key elements: (1) the order structure, which
can affect the choices of positive ideal-points and negative ideal-points;
(2) the distance/similarity measure, which is closely related to the relative
closeness degrees and determines the accuracy and rationality of decision-making.
Therefore, essentially speaking, all improvements on IF TOPSIS method provide certain improvements
on the order structure and the distance/similarity measures.

 To rank IFVs, score function is a useful tool. Xu and Yager~\cite{XY2006}
 introduced the linear order `$\leq_{_{\text{XY}}}$' for IFVs by applying
 a score function and an accuracy function. According to the TOPSIS
 idea~\cite{HY1981}, Zhang and Xu~\cite{ZX2012} proposed another linear order
 `$\leq_{_{\text{ZX}}}$' for IFVs by applying a similarity function and an
 accuracy function. Recently, Xing et al.~\cite{XXL2018} defined a linear
 order for IFVs by using a score function expressed by the Euclidean distance
 from the maximum point $\langle 1, 0\rangle$. { Bustince et al.~\cite{BFKM2013}
 suggested a general construction of linear orders for intervals contained in $[0, 1]$
 by means of aggregation functions. Based on this, De~Miguel et al.~\cite{DeBFIKM2016,
 DeBPBDaBMO2016} developed a general method for constructing linear orders
 between pairs of intervals based on aggregation functions, which was
 successfully applied to construct linear orders for interval-valued IFSs.}
 Wang et al.~\cite{WZL2006} classified the existing score functions of
 IFVs into two types, one type consists of score functions without abstention group
 influence \cite{CT1994,Kh2009,LYX2007,LW2007,WL2011,YE2007}, and the other consists of score
 functions with abstention group influence \cite{WLZ2012}. Zeng et al.~\cite{ZCK2019}
 illustrated that these existing score functions, the existing accuracy functions in
 \cite{HC2000,Ye2010}, and the measure methods in \cite{G2013,SK2009,ZX2012} for
 ranking IFVs have a drawback that they may cause some unreasonable raking results
 or they cannot distinguish some different IFVs. To overcome this drawback, Zeng et al.~\cite{ZCK2019} proposed
 a new score function $S_{_{CK}}$ for IFVs, which was a monotonically increasing
 injective with Atanassov's partial order $\subset$ (see \cite[Theorems~3.1 and 3.2]{ZCK2019}).
 However, we constructed an example to show that this inference in \cite{ZCK2019} does not hold
 and proved that such a score function does not exist, i.e., there is no any continuous
 injective from the space of all IFVs to $\mathbb{R}$ that is increasing with Atanassov's
 partial order $\subset$. This means that the score function $S_{_{CK}}$ has the same drawback.
 This trouble mainly arises from the fact that the space of all IFVs has a two-dimensional structure,
 which is not homeomorphic to any closed interval on the real line $\mathbb{R}$. In fact, we proved
 that the application of a single monotonous and continuous function does not distinguish all IFVs.

 Being a pair of dual concepts, the normalized distance measures and similarity measures
 have been widely studied. Similarly to the axiomatic definition of
 similarity measure for fuzzy sets \cite{Liu1992,Wang1982}, Li and Cheng~\cite{LC2002} gave an
 axiomatic definition of similarity measures for IFSs by using normalization (S1), symmetry (S3), and
 compatibility with Atanassov's partial order (S4), i.e., the condition $I_1\subset I_2\subset I_3$
 implies that the similarity measure between $I_1$ and $I_3$ is smaller than that between $I_1$ and $I_2$
 and between $I_2$ and $I_3$. Xu~\cite{Xu2007a} introduced some new IF similarity measures and applied
 them to some practical MADM problems. Xu and Chen~\cite{XC2008} presented a comprehensive overview
 of distance and similarity measures of IFSs and proposed some continuous distance and similarity
 measures for IFSs using the weighted Hamming distance, weighted Euclidean distance, and
 weighted Hausdorff distance. Iancu~\cite{Ia2014} defined some IF similarity measures using Frank
 t-norms $T_{\gamma }^{\mathbf{F}}$. Szmidt and Kacprzyk~\cite{SK2000} pointed out that the third
 parameter (indeterminacy degree) should be considered when calculating distances for IFSs.
 Because of the duality between distance and similarity measures, various three-dimensional
 IF distance and similarity measures including indeterminacy degrees were introduced in
 \cite{SMLXC2018,Sz2014,XC2012,YC2012}. However, Atanassov's partial order
 only reveals the magnitudes of membership degrees and non-membership degrees between two
 IFSs, and thus many three-dimensional IF similarity measures considering indeterminacy
 degrees might not meet the axiomatic condition (S4). In fact, we constructed three examples
 to show that Euclidean similarity measure {in~\cite{Sz2014,Li2014,XC2012}}, Minkowski similarity
 measure in \cite{Xu2007a,XC2012,Sz2014,Li2014}, and a modified similarity measure
 { in} \cite{XC2008} based on the idea of the above TOPSIS does not satisfy the axiomatic
 condition (S4). In particular, it should be pointed out that some existing IF distance and
 similarity measures are unreasonable for dealing with some practical decision-making
 problems. For example, Mitchell~\cite{Mi2003} showed that Li and Cheng's similarity measure
 \cite{LC2002} may lead to counter-intuitive situations in some cases. Chen et al.~\cite{CCL2016}
 showed some counterexamples to illustrate that the similarity measures in
 \cite{CC2015,HY2004,LS2003,Mi2003,Ye2010,ZY2013}
 may produce unreasonable results in some cases. As noted above, the application of a
 single monotonous and continuous function does not distinguish all IFVs. This means
 that all continuous IF similarity measures are unable to distinguish every pair
 of different IFSs. For example, when we apply a continuous IF similarity measure for pattern
 recognition, we always encounter the case that we cannot determine the classification
 result. Therefore, the comparative analysis on the indistinguishability is meaningless
 (\cite[Tables~2--5]{BA2014}, \cite[Tables~1--2]{CC2015}, \cite[Tables~1, 2, 5, 6]{CCL2016},
 \cite[Table~2]{SMLXC2018}), because all continuous IF similarity measures will also
 encounter the same indistinguishability problem. 
On the other hand, all existing distance and similarity measures are only admissible
with Atanassov's partial order, and thus decision-making results obtained by these
distance or similarity measures can only guarantee the monotonicity with the partial
order. Therefore, to obtain monotonous decision-making methods, we need to develop
new admissible distance and similarity measures with linear orders.

 During the past decade, various generalized IF TOPSIS methods have been developed.
 For example, Boran et al.~\cite{BGKA2009} first extended the TOPSIS method to IF
 group decision-making with IFV weights. Then, some similar IF TOPSIS methods with
 IFV/linguistic weights were developed \cite{AAB2021,Li2014,MDJAR2019} and widely
 applied to practical decision-making problems~\cite{AAB2021,BBM2012,BG2016,RYU2020,ZZY2020}.
 Our Examples~\ref{Exm-1} and \ref{Exm-2} below in this paper show that the TOPSIS methods in
 \cite{BGKA2009,BBM2012,Li2014,MDJAR2019,RYU2020} may yield unreasonable results
 even when dealing with the simplest decision-making problems. Chen et al.~\cite{CCL2016a} proposed
 a MADM method with crisp numerical weights based on the TOPSIS method and a new
 similarity measure for the IF situation. Zeng et al.~\cite{ZCK2019} pointed out that
 the MADM method in \cite{CCL2016a} cannot distinguish the alternatives in some
 special situations. Furthermore, our Example~\ref{Exm-3} below in this paper demonstrates
 that the MADM method in \cite{CCL2016a} is not monotonous with the linear order `$\leq_{_{\mathrm{XY}}}$'
 of Xu and Yager~\cite{XY2006}. Based on a new distance measure, Shen et al.~\cite{SMLXC2018}
 developed an extended IF TOPSIS method and applied it to credit risk evaluation.
 {Here, only some works closely related to this paper are mentioned. For more results
 on TOPSIS under the IF or interval-valued IF framework, interested readers are referred
 to \cite{Wo2016,Yue2014,LBF2017,WC2017,Ya2018,ZCF2020c,SRAG2021}.}

Inspired by the above discussions, this paper establishes
a monotonous IF TOPSIS under three popular linear orders, `$\leq_{_{\text{XY}}}$'
in \cite{XY2006}, `$\leq_{_{\text{ZX}}}$' in \cite{ZX2012}, and `$\leq_{_{A,B}}$'
in (\cite{DeBFIKM2016,DeBPBDaBMO2016,BFKM2013}). More precisely,
the main contributions of this paper are as follows:

(1) We construct some counterexamples to illustrate that Euclidean similarity measure
(\cite{Li2014,Sz2014}), Minkowski similarity measure (\cite{Xu2007a}), and
modified Euclidean similarity measure (\cite{XC2008}) do not satisfy the
axiomatic definition of IF similarity measures
(see Examples~\ref{Exm-No-Metric-1}--\ref{Exm-No-Metric-3}).

(2) We prove that there is no any continuous and injective function from the space
of all IFVs to the real line $\mathbb{R}$ that is increasing with Atanassov's partial
order `$\subset$'. This indicates the nonexistence of continuous similarity
measure distinguishing between every pair of IFVs. Therefore, the comparative
analysis on the indistinguishability for IF similarity measures is meaningless,
but unfortunately this is a common problem for all
continuous IF similarity measures.

(3) We construct three simple examples to show that some classical IF TOPSIS methods
in~\cite{Li2014,CCL2016a,BBM2012,BGKA2009,MDJAR2019,RYU2020} are not monotonous
with Atanassov's partial order `$\subset$' or the linear order `$\leq_{_{\text{XY}}}$'
(see Examples~\ref{Exm-1}--\ref{Exm-3}), which may yield counter-intuitive
results. {To overcome this limitation, by proposing three new admissible distances
with the linear order `$\leq_{_{\text{XY}}}$' or `$\leq_{_{\textrm{ZX}}}$' or
`$\leq_{_{A,B}}$', we develop a novel IF TOPSIS method and prove that it
is monotonically increasing with these three linear orders.}

(4) We provide two practical examples and a comparative analysis with other
MADM methods to illustrate the efficiency of our proposed TOPSIS method.

{The paper is organized as follows: Section~\ref{Sec-2} presents some basic definitions
related to the IFSs, including IFSs, orders for IFSs, and IF distance and similarity measure.
Section~\ref{Sec-3} provides some examples to illustrate that classical similarity measures
in \cite{Sz2014,Xu2007a,XC2008} based on the normalized Euclidean distance and normalized
Minkowski distance do not meet the axiomatic definition of IF similarity measures.
Section~\ref{Sec-4} proves that the application of a single monotonous and continuous
function does not distinguish all IFVs. Section~\ref{Sec-5} applies three examples to
demonstrate that the IF TOPSIS methods in \cite{Li2014,CCL2016a,BBM2012,BGKA2009,MDJAR2019,RYU2020}
are not monotonous with Atanassov's partial order `$\subset$' or the linear order
`$\leq_{_{\text{XY}}}$' of Xu and Yager~\cite{XY2006}. To overcome this limitation,
by constructing three new admissible distances with the linear order `$\leq_{_{\text{XY}}}$',
`$\leq_{_{\textrm{ZX}}}$' or `$\leq_{_{A,B}}$', Section~\ref{Sec-6} develops a new IF
TOPSIS method and proves that it is monotonically increasing with these three linear orders.
Section~\ref{Sec-7} presents two practical examples to demonstrate the efficiency
of the proposed TOPSIS method developed in Section~\ref{Sec-6}. Section~\ref{Sec-8}
concludes the paper with a future research outlook.}

\section{Preliminaries}\label{Sec-2}

\subsection{Intuitionistic fuzzy set (IFS)}

\begin{definition}[{\textrm{\protect\cite[Definition~1.1]{Ata1999}}}]
Let $X$ be the universe of discourse. An \textit{intuitionistic fuzzy set}~(IFS)
$I$ in $X$ is defined as an object in the following form:
\begin{equation*}
I=\left\{\langle x; \mu_{_{I}}(x), \nu_{_{I}}(x)\rangle \mid x\in X\right\},
\end{equation*}
where the functions
$\mu_{_{I}}: X \rightarrow [0,1]$
and
$\nu_{_{I}}: X \rightarrow [0,1]$
define the \textit{degree of membership} and the
\textit{degree of non-membership} of an element $x \in X$ to the set $I$,
respectively, and for every $x\in X$, $\mu_{_{I}}(x)+\nu_{_{I}}(x)\leq 1.$
\end{definition}

Let $\mathrm{IFS}(X)$ denote the set of all IFSs in $X$.
For $I\in \mathrm{IFS}(X)$, the \textit{indeterminacy degree} $\pi_{_{I}}
(x)$ of an element $x$ belonging to $I$ is defined by $\pi_{_I}(x)=1-\mu_{_I}(x)-\nu_{_I}(x)$.
In \cite{Xu2007,XC2012}, the pair $\langle\mu_{_I}(x), \nu_{_I}(x)\rangle$ is
called an \textit{ intuitionistic fuzzy value} (IFV) or an \textit{intuitionistic fuzzy number} (IFN).
For convenience, use $\alpha=\langle \mu_{\alpha}, \nu_{\alpha}\rangle$ to represent an IFV $\alpha$,
which satisfies $\mu_{\alpha}\in [0, 1]$, $\nu_{\alpha}\in [0, 1]$, and $0\leq \mu_{\alpha}
+\nu_{\alpha}\leq 1$. Additionally, $s(\alpha) =\mu_{\alpha}-\nu_{\alpha}$ and $h(\alpha)=\mu_{\alpha}+\nu_{\alpha}$
are called the \textit{score degree} and the \textit{accuracy degree} of $\alpha$, respectively.
Let $\tilde{\mathbb{I}}$ denote the set of all IFVs, i.e.,
${\tilde{\mathbb{I}}}=\{\langle \mu, \nu \rangle\in [0, 1]^{2} \mid \mu+\nu \leq 1\}$.

Motivated by the basic operations on IFSs, Xu et al.~\cite{XC2012,XY2006}
introduced the following basic operational laws for IFVs.

\begin{definition}[{\textrm{\protect\cite[Definition~1.2.2]{XC2012}}}]
\label{Def-Int-Operations}
Let $\alpha=\langle\mu_{\alpha}, \nu_{\alpha}\rangle$,
$\beta=\langle\mu_{\beta}, \nu_{\beta}\rangle\in \tilde{\mathbb{I}}$. Define
\begin{enumerate}[(i)]
\item $\alpha^{\complement}=\langle\nu_{\alpha},\mu_{\alpha}\rangle$.
\item $\alpha\cap\beta=\langle\min\{\mu_{\alpha}, \mu_{\beta}\}, \max\{\nu_{\alpha},\nu_{\beta}\}\rangle$.
\item $\alpha\cup\beta=\langle\max\{\mu_{\alpha},\mu_{\beta}\}, \min\{\nu_{\alpha},\nu_{\beta}\}\rangle$.
\item $\alpha\oplus\beta=\langle\mu_{\alpha}+\mu_{\beta}-\mu_{\alpha}\mu_{\beta}, \nu_{\alpha}\nu_{\beta}\rangle$.
\item $\alpha\otimes\beta=\langle\mu_{\alpha}\mu_{\beta}, \nu_{\alpha}+\nu_{\beta}-\nu_{\alpha}\nu_{\beta}\rangle$.
\item $\lambda\alpha=\langle 1-(1-\mu_{\alpha})^{\lambda}, (\nu_{\alpha})^{\lambda}\rangle$, $\lambda >0 $.
\item $\alpha^{\lambda}=\langle (\mu_{\alpha})^{\lambda}, 1-(1-\nu_{\alpha})^{\lambda}\rangle$, $\lambda >0$.
\end{enumerate}
\end{definition}

\subsection{Orders for IFSs}

Atanassov's order `$\subset$' \cite{Ata1999}, defined by that $\alpha\subset \beta$ if and only if
$\alpha\cap \beta=\alpha$, is a partial order on $\tilde{\mathbb{I}}$.
{\begin{definition}[{\textrm{\protect\cite[Definition 4.1]{DeBFIKM2016}}}]
\label{de-order(ZX)}
An order $\leq$ on $\tilde{\mathbb{I}}$ is said to be an \textit{IF-admissible order}
if it is a linear order and refines Atanassov's order $\subset$.
\end{definition}}

To compare any two IFVs, Xu and Yager~\cite{XY2006} introduced the following
linear order `$\leq_{_{\text{XY}}}$' (see also \cite[Definition~3.1]{Xu2007}
and \cite[Definition~1.1.3]{XC2012}):

\begin{definition}
[{\textrm{\protect\cite[Definition~1]{XY2006}}}]
\label{de-order(Xu)}
Let $\alpha_{1}$ and $\alpha_{2}$ be two IFVs.
\begin{itemize}
  \item If $s(\alpha_{1})<s(\alpha_{2})$,
then $\alpha_{1}$ is smaller than $\alpha_{2}$,
denoted by $\alpha_{1}<_{_{\text{XY}}}\alpha_{2}$.

  \item If $s(\alpha_{1})=s(\alpha_{2})$, then

\begin{itemize}
  \item if $h(\alpha_{1})=h(\alpha_{2})$, then $\alpha_{1}=\alpha_{2}$;
  \item if $h(\alpha_{1})<h(\alpha_{2})$, then $\alpha_{1}<_{_{\text{XY}}}\alpha_{2}$.
\end{itemize}
\end{itemize}
If $\alpha_{1}<_{_{\text{XY}}}\alpha_{2}$ or $\alpha_{1}=\alpha_{2}$, then denote it by $\alpha_{1}\leq_{_{\text{XY}}}\alpha_{2}$.
\end{definition}

Alongside Xu and Yager's order `$\leq_{_{\text{XY}}}$' in Definition \ref{de-order(Xu)},
Szmidt and Kacprzyk~\cite{SK2009} proposed another comparison
function, $\rho(\alpha)=\frac{1}{2}(1+\pi(\alpha))(1-\mu(\alpha))$ for IFVs,
which is a partial order. However, it sometimes cannot distinguish between two IFVs. Although Xu's method
\cite{Xu2007} constructs a linear order for ranking any two IFVs, its procedure has
the following disadvantages:
(1) It may result in that the less we know, the better the IFV, which is not reasonable.
(2) It is sensitive to a slight change of the parameters.
(3) It is not preserved under multiplication by a scalar, namely,
$\alpha\leq_{_{\text{XY}}}\beta$ might not imply
$\lambda\alpha<_{_{\text{XY}}}\lambda\beta$, where $\lambda$ is a scalar
(see {\textrm{\protect\cite[Example~1]{BBGMP2011}}}).
To overcome such shortcomings of the above two ranking methods, Zhang and Xu~\cite{ZX2012}
improved Szmidt and Kacprzyk's method \cite{SK2009}, according to Hwang and Yoon's idea
\cite{HY1981} and technique for preference order by similarity to an ideal point.
They also introduced a \textit{similarity function} $L(\alpha)$,
called the \textit{$L$-value} in \cite{ZX2012}, for any IFV
$\alpha=\langle\mu_{\alpha}, \nu_{\alpha}\rangle$, as follows:
\begin{equation}
L(\alpha)=\frac{1-\nu_{\alpha}}{(1-\mu_{\alpha})+(1-\nu_{\alpha})}
=\frac{1-\nu_{\alpha}}{1+\pi_{\alpha}}.
\end{equation}
In particular, if $\nu_{\alpha}<1$, then
$L(\alpha)=\frac{1}{\tfrac{1-\mu_{\alpha}}{1-\nu_{\alpha}}+1}.$
Furthermore, Zhang and Xu~\cite{ZX2012} introduced the following order `$\leq_{_{\textrm{ZX}}}$'
for IFVs by applying the similarity function $L(\_)$.

\begin{definition}[{\textrm{\protect\cite{ZX2012}}}]
\label{de-order(ZX)}
Let $\alpha_{1}$ and $\alpha_{2}$ be two IFVs.
\begin{itemize}
  \item If $L(\alpha_{1})<L(\alpha_{2})$, then $\alpha_{1}<_{_{\textrm{ZX}}}\alpha_{2}$;
  \item If $L(\alpha_{1})=L(\alpha_{2})$, then
  \begin{itemize}
    \item if $h(\alpha_{1})=h(\alpha_{2})$, then $\alpha_{1}=\alpha_{2}$;
    \item if $h(\alpha_{1})<h(\alpha_{2})$, then $\alpha_{1}<_{_{\textrm{ZX}}}\alpha_{2}$.
  \end{itemize}
\end{itemize}
If $\alpha_{1}<_{_{\textrm{ZX}}}\alpha_{2}$ or $\alpha_{1}=\alpha_{2}$,
then denote it by $\alpha_{1}\leq_{_{\textrm{ZX}}}\alpha_{2}$.
\end{definition}

{

\begin{definition}
[{\textrm{\protect\cite[Definition~1.1]{GMMP2009}}}]
\label{de-Agg-Func}
An \textit{aggregation function} in $[0, 1]^{n}$ is a function
$A: [0, 1]^{n}\rightarrow [0, 1]$ that
\begin{enumerate}[{\rm (i)}]
  \item is nondecreasing in each variable;
  \item fulfills the boundary conditions
  $$
  A(0, \ldots, 0)=0 \text{ and }
  A(1, \ldots, 1)=1.
  $$
\end{enumerate}
\end{definition}

By an equivalent transformation between intervals and IFVs and \cite[Proposition~3.2]{BFKM2013},
the following general construction of linear orders is proposed for IFVs.
\begin{proposition}
\label{Prop-Linear-Order}
Let $A$, $B:[0, 1]^{2}\rightarrow [0, 1]$ be two continuous aggregation functions satisfying
that, for $(x_1, y_1)$, $(x_2, y_2)\in [0, 1]^{2}$, the identities $A(x_1, y_1)= A(x_2, y_2)$ and
$B(x_1, y_1)= B(x_2, y_2)$ can hold only if $x_1=x_2$ and $y_1=y_2$. Define the order
$\leq_{_{A,B}}$ on $\tilde{\mathbb{I}}$ as follows:
$\langle \mu_1, \nu_1\rangle \leq_{_{A,B}}\langle \mu_2, \nu_2\rangle$
if and only if $A(\mu_1, 1-\nu_1)< A(\mu_2, 1-\nu_2)$ or $(A(\mu_1, 1-\nu_1)=
A(\mu_2, 1-\nu_2)$ and $B(\mu_1, 1-\nu_1)\leq B(\mu_2, 1-\nu_2))$.
Then, $\leq_{_{A,B}}$ is an admissible order on $\tilde{\mathbb{I}}$.
\end{proposition}

For $\gamma\in [0, 1]$, consider an aggregation function $K_{\gamma}:
[0, 1]^2\rightarrow [0, 1]$ defined by $K_{\gamma}(x, y)=x+\gamma (y-x).$
For $\gamma_1$, $\gamma_2\in [0, 1]$ with $\gamma_1\neq \gamma_2$,
according to Proposition~\ref{Prop-Linear-Order},
it follows that the order $\leq_{_{\gamma_1, \gamma_2}}$ on $\tilde{\mathbb{I}}$
defined by $\alpha \leq_{_{\gamma_1, \gamma_2}} \beta$ if and only if
$K_{\gamma_1}(\mu_{\alpha}, 1-\nu_{\alpha})<K_{\gamma_1}(\mu_{\beta}, 1-\nu_{\beta})$
or ($K_{\gamma_1}(\mu_{\alpha}, 1-\nu_{\alpha})=K_{\gamma_1}(\mu_{\beta}, 1-\nu_{\beta})$
and $K_{\gamma_2}(\mu_{\alpha}, 1-\nu_{\alpha})\leq K_{\gamma_2}(\mu_{\beta}, 1-\nu_{\beta})$)
is an admissible order on $\tilde{\mathbb{I}}$.
}

\subsection{IF distance and similarity measure}

Li and Cheng~\cite{LC2002} introduced an axiomatic definition of similarity measure for IFSs, which
was then improved by Mitchell~\cite{Mit2003} as follows. More results on the similarity measure
can be found in \cite{Sz2014,BurB1996}.

\begin{definition}[{\textrm{\protect\cite{Mit2003}}}]
\label{Def-Li-Cheng}
Let $X$ be the universe of discourse. A mapping $\mathbf{S}: \mathrm{IFS}(X)\times \mathrm{IFS}(X)
\rightarrow [0, 1]$ is called an \textit{admissible similarity measure
with the order $\subset$} on $\mathrm{IFS}(X)$ if it satisfies the following conditions:
for any $I_1$, $I_2$, $I_3\in \mathrm{IFS}(X)$,
\begin{enumerate}
  \item[(S1)] $0\leq \mathbf{S}(I_1, I_2)\leq 1$.
  \item[(S2)] $\mathbf{S}(I_1, I_2)=1$ if and only if $I_1=I_2$.
  \item[(S3)] $\mathbf{S}(I_1, I_2)=\mathbf{S}(I_2, I_1)$.
  \item[(S4)] If $I_1\subset I_2\subset I_3$, then $\mathbf{S}(I_1, I_3)\leq \mathbf{S}(I_1, I_2)$
  and $\mathbf{S}(I_1, I_3)\leq \mathbf{S}(I_2, I_3)$.
\end{enumerate}
\end{definition}

\begin{remark}
The admissible similarity measure with the order $\subset$ was also called \textit{similarity measure}
by Hung and Yang~\cite{HY2008} and Szmidt~\cite{Sz2014}. When no ambiguity is possible, we simply
call it similarity measure.
\end{remark}

\begin{definition}
\label{Def-Wu-1}
Let $X$ be the universe of discourse and $I_1$, $I_2\in \mathrm{IFS}(X)$.
If $\langle \mu_{_{I_1}}(x), \nu_{_{I_1}}(x)\rangle \leq_{_{\text{XY}}}
\langle \mu_{_{I_2}}(x), \nu_{_{I_2}}(x)\rangle$ holds for all $x\in X$, then
we say that $I_{1}$ is \textit{smaller} than or equal to $I_2$ under the
linear order $\leq_{_{\text{XY}}}$, denoted by $I_{1} \leq_{_{\text{XY}}} I_2$.
\end{definition}

Based on Definition~\ref{Def-Wu-1},
we introduce the improved similarity measure definition for IFSs below:

\begin{definition}
\label{Def-Wu-Sim}
Let $X$ be the universe of discourse. A mapping $\mathbf{S}: \mathrm{IFS}(X)\times \mathrm{IFS}(X)
\rightarrow [0, 1]$ is called an {\it admissible similarity measure
with the order $\leq_{_{\text{XY}}}$} on $\mathrm{IFS}(X)$ if it satisfies the
conditions (S1)--(S3) in Definition~\ref{Def-Li-Cheng}, and the following one (S4$^{\prime}$):
\begin{enumerate}
  \item[(S4$^{\prime}$)] For any $I_1$, $I_2$, $I_3\in \mathrm{IFS}(X)$,
  if $I_1\leq_{_{\text{XY}}} I_2\leq_{_{\text{XY}}} I_3$, then $\mathbf{S}(I_1, I_3)\leq \mathbf{S}(I_1, I_2)$
  and $\mathbf{S}(I_1, I_3)\leq \mathbf{S}(I_2, I_3)$.
\end{enumerate}
\end{definition}

Now, we recall some classical distances and similarity measures for IFSs.

The \textit{normalized Hamming distance} in \cite{SK2000} is:
\begin{equation}\label{Dis-Hamming}
d_{\mathrm{Ha}}(I_1, I_2)=\frac{1}{2n}\sum_{j=1}^{n}H_{j},
\end{equation}
where $H_{j}=|\mu_{I_1}(x_j)-\mu_{I_2}(x_j)|
+|\nu_{I_1}(x_j)-\nu_{I_2}(x_j)|+|\pi_{I_1}(x_j)-\pi_{I_2}(x_j)|$.

The \textit{normalized Euclidean distance} in \cite{SK2000} is:
\begin{equation}\label{Dis-Euclidean}
d_{\mathrm{Eu}}(I_1, I_2)=\sqrt{\frac{1}{2n}\sum_{j=1}^{n}E_{j}},
\end{equation}
where $E_{j}=|\mu_{I_1}(x_j)-\mu_{I_2}(x_j)|^2
+|\nu_{I_1}(x_j)-\nu_{I_2}(x_j)|^2+|\pi_{I_1}(x_j)-\pi_{I_2}(x_j)|^2$.

The \textit{normalized Minkowski distance} in \cite{SK2000,Sz2014} is:
\begin{equation}\label{Dis-Minkowski}
d_{_{\mathrm{M}}}^{(\alpha)}(I_1, I_2)=\sqrt[\alpha]{\frac{1}{2n}
\sum_{j=1}^{n}M_{j}},
\end{equation}
where $M_{j}=|\mu_{I_1}(x_j)-\mu_{I_2}(x_j)|^{\alpha}
+|\nu_{I_1}(x_j)-\nu_{I_2}(x_j)|^{\alpha}+|\pi_{I_1}(x_j)-\pi_{I_2}(x_j)|^{\alpha}$
and $\alpha\geq 1$.

By using the normalized Hamming distance, Szmidt and Kacprzyk~\cite{SK2004} introduced the
following similarity measure $\mathbf{S}_{\mathrm{SK}}^{1}$ for two IFSs $I_1$ and $I_2$:
\begin{equation}
\label{SK-similarity}
\mathbf{S}_{_{\mathrm{SK}}}^{1}(I_1, I_2)=
\frac{d_{\mathrm{Ha}}(I_1, I_2)}{d_{\mathrm{Ha}}(I_1, I_2^{\complement})},
\end{equation}
where $I_2^{\complement}$ is the complement of $I_2$.
If we replace the normalized Hamming distance with the normalized Euclidean distance,
we can obtain the following ``similarity measure" $\mathbf{S}_{\mathrm{SK}}^{1}$ for two IFSs $I_1$ and $I_2$:
\begin{equation}
\label{SK-similarity-2}
\mathbf{S}_{_{\mathrm{SK}}}^{2}(I_1, I_2)=
\frac{d_{\mathrm{Eu}}(I_1, I_2)}{d_{\mathrm{Eu}}(I_1, I_2^{\complement})}.
\end{equation}
Clearly, both $\mathbf{S}_{\mathrm{SK}}^{1}$ and $\mathbf{S}_{\mathrm{SK}}^{2}$
are not similarity measures, because their values may exceed $1$. Then, the following
similarity measures were introduced by Szmidt~\cite{Sz2014} using the normalized Hamming
and Euclidean distances:
\begin{equation}
\label{S-similarity}
\mathbf{Sim}_{_{\text{H}}}^{1}(I_1, I_2)=1-d_{\mathrm{Ha}}(I_1, I_2),
\end{equation}
and
\begin{equation}
\label{S-similarity-2}
\mathbf{Sim}_{_{\text{E}}}^{2}(I_1, I_2)=1-d_{\mathrm{Eu}}(I_1, I_2).
\end{equation}

Based on the normalized Minkowski distance $d_{_{\mathrm{M}}}^{(\alpha)}$,
Xu~\cite{Xu2007a} introduced the following \textit{Minkowski similarity measure}
(see also \cite{Li2014,XC2008}):
for $\alpha>0$,
\begin{equation}
\label{Xu-similarity-3}
\mathbf{S}_{_{\text{Xu}}}^{(\alpha)}(I_1, I_2)=1-d_{_{\mathrm{M}}}^{(\alpha)}(I_1, I_2).
\end{equation}
Clearly, $\mathbf{S}_{_{\text{Xu}}}^{(1)}=\mathbf{Sim}_{_{\text{H}}}^{1}$
and $\mathbf{S}_{_{\text{Xu}}}^{(2)}=\mathbf{Sim}_{_{\text{E}}}^{2}$.

Motivated by the idea of the TOPSIS of Hwang and Yoon~\cite{HY1981},
Xu and Chen~\cite{XC2008} modified Eqs.~\eqref{SK-similarity} and \eqref{SK-similarity-2} as follows:
\begin{equation}
\label{Xu-similarity}
\mathbf{S}_{_{\text{XC}}}^{1}(I_1, I_2)=
\frac{d_{\mathrm{Ha}}(I_1, I_2^{\complement})}{d_{\mathrm{Ha}}(I_1, I_2)+
d_{\mathrm{Ha}}(I_1, I_2^{\complement})},
\end{equation}
\begin{equation}
\label{Xu-similarity-2}
\mathbf{S}_{_{\text{XC}}}^{2}(I_1, I_2)=
\frac{d_{\mathrm{Eu}}(I_1, I_2^{\complement})}{d_{\mathrm{Eu}}(I_1, I_2)+
d_{\mathrm{Eu}}(I_1, I_2^{\complement})}.
\end{equation}


\section{The drawbacks of some existing similarity measures}\label{Sec-3}

This section illustrates that the similarity measures defined
by Eqs.~\eqref{S-similarity-2}, \eqref{Xu-similarity-3}, and
\eqref{Xu-similarity-2} do not meet the property (S4) in
the axiomatic definition of intuitionistic fuzzy similarity measures.

\begin{example}
\label{Exm-No-Metric-1}
Let the universe of discourse $X=\{x_1\}$, and
$I_1=\left\{\frac{\langle 0, 1\rangle}{x_1}\right\}$,
$I_2=\left\{\frac{\langle 0.1, 0\rangle}{x_1}\right\}$, and
$I_3=\left\{\frac{\langle 0.4, 0\rangle}{x_1}\right\}$.
Clearly, $I_1\subset I_2\subset I_3$. By direct calculation,
we have
\begin{align*}
\mathbf{Sim}_{_{\text{E}}}^{2}(I_1, I_2)
=& 1-\sqrt{\frac{|0-0.1|^2+|1-0|^2+|0-0.9|^2}{2}}\\
=& 1-\sqrt{0.91},
\end{align*}
and
\begin{align*}
\mathbf{Sim}_{_{\text{E}}}^{2}(I_1, I_3)
=& 1-\sqrt{\frac{|0-0.4|^2+|1-0|^2+|0-0.6|^2}{2}}\\
=& 1-\sqrt{0.76},
\end{align*}
and thus $\mathbf{Sim}_{_{\text{E}}}^{2}(I_1, I_2)<\mathbf{Sim}_{_{\text{E}}}^{2}(I_1, I_3)$.
This, together with $I_1\subset I_2\subset I_3$, implies that
the formula $\mathbf{Sim}_{_{\text{H}}}^{1}$ defined by Eq.~\eqref{S-similarity-2}
is not a similarity measures on IFSs.
\end{example}

\begin{example}
Let the universe of discourse $X=\{x_1\}$ and $I_1=\left\{\frac{\langle 0, 1\rangle}{x_1}\right\}$,
$I_2=\left\{\frac{\langle \mu_2, 0\rangle}{x_1}\right\}$, and $I_3=\left\{\frac{\langle \mu_3, 0\rangle}{x_1}\right\}$
be three IFSs on $X$ such that $0<\mu_2<\mu_3<0.5$. Clearly, $I_1\subset I_2\subset I_3$. Fix any $\alpha>1$.
By direct calculation, we have
$$
\mathbf{S}_{_{\text{M}}}^{(\alpha)}(I_1, I_2)
=1-\sqrt[\alpha]{\frac{(\mu_2)^{\alpha}+1+(1-\mu_2)^{\alpha}}{2}},
$$
and
$$
\mathbf{S}_{_{\text{M}}}^{(\alpha)}(I_1, I_3)
=1-\sqrt[\alpha]{\frac{(\mu_3)^{\alpha}+1+(1-\mu_3)^{\alpha}}{2}}.
$$
Let $\Gamma(x)=1-\sqrt[\alpha]{\frac{x^{\alpha}+1+(1-x)^{\alpha}}{2}}$ ($x\in (0, 0.5)$).
Noting that $\alpha>1$ and $x\in (0, 0.5)$, by direct calculation, we get $\Gamma^{\prime}(x)=-\frac{1}{2}(\frac{x^{\alpha}+1+(1-x)^{\alpha}}{2})
^{\frac{1}{\alpha}-1}\cdot (x^{\alpha-1}-(1-x)^{\alpha-1})>0$, and thus
the function $\Gamma$ is strictly increasing on $(0, 0.5)$. This, together
with $0<\mu_2<\mu_3<0.5$, implies that $\mathbf{S}_{_{\text{M}}}^{(\alpha)}(I_1, I_2)
=\Gamma(\mu_2)<\Gamma(\mu_3)=\mathbf{S}_{_{\text{M}}}^{(\alpha)}(I_1, I_3)$.
Therefore, the formula $\mathbf{S}_{_{\text{M}}}^{(\alpha)}$ defined by
Eq.~\eqref{Xu-similarity-3} is not a similarity measures on IFSs for any $\alpha>1$.
\end{example}

\begin{example}
\label{Exm-No-Metric-3}
Let the universe of discourse $X=\{x_1\}$, and $I_1=\left\{\frac{\langle 0, 1\rangle}{x_1}\right\}$,
$I_2=\left\{\frac{\langle 0.9, 0.01\rangle}{x_1}\right\}$, and
$I_3=\left\{\frac{\langle 0.901, 0.007\rangle}{x_1}\right\}$.
Clearly, $I_1\subset I_2\subset I_3$. By direct calculation, we have
{\small\begin{align*}
& \mathbf{S}_{_{\text{XC}}}^{2}(I_1, I_2)\\
=&
\frac{\sqrt{\frac{|0-0.01|^2+|1-0.9|^2+|0-0.09|^2}{2}}}
{\sqrt{\frac{|0-0.9|^2+|1-0.01|^2+|0-0.09|^2}{2}}+
\sqrt{\frac{|0-0.01|^2+|1-0.9|^2+|0-0.09|^2}{2}}}\\
\approx & 0.09141,
\end{align*}}
and
{\small \begin{align*}
& \mathbf{S}_{_{\text{XC}}}^{2}(I_1, I_3)\\
=&
\frac{\sqrt{\frac{|0.007|^2+|1-0.901|^2+|0.092|^2}{2}}}
{\sqrt{\frac{|0.901|^2+|1-0.007|^2+|0.092|^2}{2}}
+\sqrt{\frac{|0.007|^2+|1-0.901|^2+|0.092|^2}{2}}}\\
\approx & 0.09148,
\end{align*}}
and thus $\mathbf{S}_{_{\text{XC}}}^{2}(I_1, I_2)<\mathbf{S}_{_{\text{XC}}}^{2}(I_1, I_3)$.
This, together with $I_1\subset I_2\subset I_3$, implies that
the formula $\mathbf{S}_{_{\text{XC}}}^{2}$ defined by Eq.~\eqref{Xu-similarity-2}
is not a similarity measures on IFSs.
\end{example}

\section{A remark on score functions for IFVs}\label{Sec-4}

Zeng et al.~\cite{ZCK2019} introduced the following {\it score value} $S_{_{CK}}(\_)$ for IFV $\alpha=
\langle \mu_{\alpha}, \nu_{\alpha}\rangle$:
{\small\begin{equation}
S_{_{CK}}(\alpha)=(\mu_{\alpha}-\nu_{\alpha})-(1-\mu_{\alpha}-\nu_{\alpha})\times
\frac{\log_{2} (2-\mu_{\alpha}-\nu_{\alpha})}{100}.
\end{equation}}
Then, they proved the following basic properties for the score function $S_{_{CK}}(\_)$.
\begin{theorem}[{\textrm{\protect\cite[Theorem~3.1]{ZCK2019}}}]
\label{ZCK-Thm}
Assume that $\alpha$ and $\beta$ are two IFVs.
If $\alpha\neq \beta$, then $S_{_{CK}}(\alpha)\neq S_{_{CK}}(\beta)$.
\end{theorem}

\begin{theorem}[{\textrm{\protect\cite[Theorem~3.2]{ZCK2019}}}]
\label{ZCK-Mono-Thm}
Assume that $\alpha$ and $\beta$ are two IFVs.
If $\alpha\supset \beta$, then $S_{_{CK}}(\alpha)> S_{_{CK}}(\beta)$.
\end{theorem}

However, the following example shows that Theorem~\ref{ZCK-Thm} does not hold.

\begin{example}
Choose $\alpha=\langle 0, 0\rangle$ and $\beta=\langle \frac{99}{200}, \frac{101}{200}\rangle$.
By direct calculation, we have
$S_{_{CK}}(\alpha)=(0-0)-(1-0-0)\times \frac{\log_{2}2}{100}=-\frac{1}{100}$
and
$S_{_{CK}}(\beta)=(\frac{99}{200}-\frac{101}{200})-(1-\frac{99}{200}
-\frac{101}{200})\times \frac{\log_{2}(2-1)}{100}
=-\frac{1}{100}.$
This implies that Theorem~\ref{ZCK-Thm} does not hold since $\alpha\neq \beta$.
\end{example}

In fact, we can prove that there is no any continuous
function from $\tilde{\mathbb{I}}$ to $\mathbb{R}$
simultaneously meeting the conditions in Theorems~\ref{ZCK-Thm} and \ref{ZCK-Mono-Thm},
which indicates that the two-dimensional structure of IFVs
is too complex to distinguish all IFVs with only a single monotonous and continuous
function, where the monotonicity is under Atanassov's order `$\subset$', and the continuity
is under the topology of subset of $\mathbb{R}^2$.

\begin{theorem}
\label{No-exist-Thm}
There is no any continuous function $f: \tilde{\mathbb{I}}\rightarrow \mathbb{R}$
satisfying the following two conditions:
\begin{enumerate}[{\rm (1)}]
  \item $f$ is injective, i.e., for any $\alpha$, $\beta
\in \tilde{\mathbb{I}}$ with $\alpha\neq \beta$, $f(\alpha)\neq f(\beta)$;

  \item $f$ is increasing under the partial order $\subset$, i.e., for any $\alpha$, $\beta\in
 \tilde{\mathbb{I}}$ with $\alpha\subset \beta$, $f(\alpha)\leq f(\beta)$.
\end{enumerate}
\end{theorem}

\begin{IEEEproof}
Suppose on the contrary that there exists a continuous function
$f: \tilde{\mathbb{I}}\rightarrow \mathbb{R}$ simultaneously satisfying
the conditions (1) and (2). 

(i) Let $\varphi(\nu)=f(\langle 0.25, \nu\rangle)$ ($\nu\in [0.25, 0.75]$). Clearly,
$\varphi$ is continuous since $f$ is continuous on $\tilde{\mathbb{I}}$.
For any $0.25\leq \nu_1\leq \nu_2\leq 0.75$, by $\langle 0.25, \nu_1\rangle
\supset \langle 0.25, \nu_2\rangle$ and condition (2), one has
$\varphi(\nu_1)=f(\langle 0.25, \nu_1\rangle)\geq f(\langle 0.25, \nu_2\rangle)
=\varphi(\nu_2)$. This, together with condition (1), implies that
$\varphi(\_)$ is strictly decreasing on $[0.25, 0.75]$. Thus,
$\varphi((0.25, 0.75])=[\varphi(0.75), \varphi(0.25))=[f(\langle 0.25, 0.75\rangle),
f(\langle 0.25, 0.25\rangle))$ by the intermediate value theorem.

  (ii) Let $\psi(\nu)=f(\langle \mu, 0.25\rangle)$ ($\mu \in [0, 0.25]$). Clearly,
$\psi$ is continuous since $f$ is continuous on $\tilde{\mathbb{I}}$.
For any $0\leq \mu_1\leq \mu_2\leq 0.25$, by $\langle \mu_1, 0.25 \rangle
\subset \langle \mu_2, 0.25 \rangle$ and condition (2), one has
$\psi(\mu_1)=f(\langle \mu_1, 0.25 \rangle)\leq f(\langle \mu_2, 0.25\rangle)
=\psi(\mu_2)$. This, together with condition (1), implies that
$\psi(\_)$ is strictly increasing on $[0, 0.25]$. Thus,
$\psi([0, 0.25))=[\psi(0), \psi(0.25))=[f(\langle 0, 0.25\rangle),
f(\langle 0.25, 0.25\rangle))$ by the intermediate value theorem.

Summing (i) and (ii), one can easily verify that
$\Lambda=\varphi((0.25, 0.75])\cap\psi([0, 0.25))
=[\max\{f(\langle 0, 0.25\rangle), f(\langle 0.25, 0.75\rangle)\},
f(\langle 0.25, 0.25\rangle))$
is a non-degenerate interval, i.e.,
$$
\max\{f(\langle 0, 0.25\rangle), f(\langle 0.25, 0.75\rangle)\}
<f(\langle 0.25, 0.25\rangle),
$$
implying that, for any $\xi\in \Lambda$, there exist $\nu\in (0.25, 0.75]$
and $\mu\in [0, 0.25)$ such that $\varphi(\nu)=\xi$ and $\psi(\mu)=\xi$,
and thus there exist $0<\mu<0.25<\nu\leq 0.75$ such that $\varphi(\nu)=f(\langle 0.25, \nu\rangle)
=\xi=f(\langle \mu, 0.25\rangle)=\psi(\mu)$. This contradicts with condition (1).
\end{IEEEproof}

Shen et al.~\cite{SMLXC2018} pointed out that many existing distance measures
cannot determine the classification results for some pattern recognition problems
(see \cite[Table~2]{SMLXC2018}), i.e., their dual similarity measures cannot
distinguish between some pair of IFVs. To overcome this drawback, they proposed a new
distance measure $d_{_{\mathrm{Sh}}}$ as follows: for $\alpha=\langle \mu_{\alpha},
\nu_{\alpha}\rangle$, $\beta=\langle \mu_{\beta},
\nu_{\beta}\rangle \in \tilde{\mathbb{I}}$,
\begin{align*}
d_{_{\mathrm{Sh}}}(\alpha, \beta)=
\sqrt{\frac{(\tilde{\mu}_{\alpha}-
\tilde{\mu}_{\beta})^2+
(\tilde{\nu}_{\alpha}-\tilde{\nu}_{\beta})^2}{2}},
\end{align*}
where $\tilde{\mu}_{\alpha}=\mu_{\alpha}(1+\frac{2}{3}\pi_{\alpha}(1+\pi_{\alpha}))$,
$\tilde{\nu}_{\alpha}=\nu_{\alpha}(1+\frac{2}{3}\pi_{\alpha}(1+\pi_{\alpha}))$,
$\tilde{\mu}_{\beta}=\mu_{\beta}(1+\frac{2}{3}\pi_{\beta}(1+\pi_{\beta}))$,
and $\tilde{\nu}_{\beta}=\nu_{\beta}(1+\frac{2}{3}\pi_{\beta}(1+\pi_{\beta}))$.

Fix $\beta\in \tilde{\mathbb{I}}$ and define $\mathcal{G}(\alpha)
=1-d_{_{\mathrm{Sh}}}(\alpha, \beta)$ for $\alpha=\langle \mu_{\alpha},
\nu_{\alpha}\rangle \in \tilde{\mathbb{I}}$. From~\cite[Theorem~1]{SMLXC2018},
it follows that, for $\alpha_1$ $\alpha_2\in \tilde{\mathbb{I}}$ with
$\alpha_1\subset \alpha_2$, one has $\mathcal{G}(\alpha_1)\leq \mathcal{G}
(\alpha_2)$, i.e., the function $\mathcal{G}$ satisfies the condition (2)
of Theorem~\ref{No-exist-Thm}. This, together with Theorem~\ref{No-exist-Thm}
and the continuity of $\mathcal{G}$, implies that $\mathcal{G}$ is not injective,
and thus there exist two different IFVs $\alpha_1$ and $\alpha_2 \in \tilde{\mathbb{I}}$
such that $\mathcal{G}(\alpha_1)=\mathcal{G}(\alpha_2)$, implying that
$d_{_{\mathrm{Sh}}}(\alpha_1, \beta)=d_{_{\mathrm{Sh}}}(\alpha_2, \beta)$.
Therefore, the distance measure $d_{_{\mathrm{Sh}}}$ has the same drawback
(see Example~\ref{Exm-SMLXC}). In fact, by Theorem~\ref{No-exist-Thm}, we
conclude that there is no any continuous distance measure that can overcome
the above drawback. Therefore, the comparative analysis in~\cite[Tables~2--5]{BA2014},
\cite[Tables~1--2]{CC2015}, \cite[Tables~1, 2, 5, 6]{CCL2016}, and \cite[Table~2]{SMLXC2018}
on the indistinguishability is meaningless.

\begin{example}
\label{Exm-SMLXC}
Let $\beta=\langle 0, 0 \rangle$ and $\alpha=\langle x, y\rangle\in \tilde{\mathbb{I}}$. Then,
\begin{align*}
&d_{_{\mathrm{Sh}}}(\alpha, \beta)\\
=&
\sqrt{\frac{[1+\frac{2}{3}(1-x-y)(2-x-y)]^{2}\times (x^2+y^2)}{2}}=0.5,
\end{align*}
i.e.,
\begin{equation}
\label{eq-Inf}
\left[1+\frac{2}{3}(1-x-y)(2-x-y)\right]^{2}\times (x^2+y^2)=0.5.
\end{equation}
Clearly, Eq.~\eqref{eq-Inf} has infinitely many solutions.
This means that there exist infinitely many IFVs, whose distances
from $\beta$ are all equal to $0.5$.
\end{example}

\section{A monotonous IF TOPSIS method with the linear orders $\leq_{_\text{XY}}$
and $\leq_{_\text{ZX}}$}\label{Sec-5}


Suppose that there are $n$ alternatives ${A}_{i}$ ($i=1, 2, \ldots, n$) evaluated with
respect to $m$ attributes $\mathscr{O}_j$ ($j=1, 2, \ldots, m$). The sets of the
alternatives and attributes are denoted by $A=\{A_1, A_2, \ldots, A_n\}$
and $\mathscr{O}=\{\mathscr{O}_1, \mathscr{O}_2, \ldots, \mathscr{O}_m\}$,
respectively. The rating (or evaluation) of each alternative $A_i\in A$ ($i=1, 2, \ldots, n$)
on each attribute $\mathscr{O}_j$ ($j=1, 2, \ldots, m$) is expressed with an IFS $r_{ij}=
\left\{\frac{\langle\mu_{ij}, \nu_{ij}\rangle}{(A_i, \mathscr{O}_j)}\right\}$,
denoted by $r_{ij}=\langle\mu_{ij}, \nu_{ij}\rangle$ for short, where
$\mu_{ij}\in [0, 1]$ and $\nu_{ij}\in [0, 1]$ are respectively the satisfaction
(or membership) degree and dissatisfaction (or non-membership) degree of the alternative
$A_{i}\in A$ on the attribute $\mathscr{O}_{j}$ satisfying the condition $0\leq \mu_{ij}+\nu_{ij}\leq 1$.
A multi-attribute decision-making (MADM) problem with IFSs is expressed in matrix form
shown in Table~\ref{Tab-IFDM}.
\begin{table}[H]	
	\centering
	\caption{IF decision matrix $R=(r_{ij})_{n\times m}$}
	\label{Tab-IFDM}
     \scalebox{0.85}{
	\begin{tabular}{cccccc}
		\toprule
		$$ & $\mathscr{O}_{1}$ &  $\mathscr{O}_{2}$ & $\ldots$ & $\mathscr{O}_{m}$ \\
		\midrule
		$A_{1}$ & $\langle \mu_{11}, \nu_{11}\rangle$ &  $\langle \mu_{12}, \nu_{12}\rangle$ & $\ldots$ & $\langle \mu_{1m}, \nu_{1m}\rangle$ \\
		$A_{2}$ & $\langle \mu_{21}, \nu_{21}\rangle$ &  $\langle \mu_{22}, \nu_{22}\rangle$ & $\ldots$ & $\langle \mu_{2m}, \nu_{2m}\rangle$ \\
        $\vdots$ & $\vdots$ &  $\vdots$ & $\ddots$ & $\vdots$ \\
		$A_{n}$ & $\langle \mu_{n1}, \nu_{n1}\rangle$ &  $\langle \mu_{n2}, \nu_{n2}\rangle$ & $\ldots$ & $\langle \mu_{nm}, \nu_{nm}\rangle$ \\
        \bottomrule
	\end{tabular}
      }
\end{table}

To follow the common sense, a good method for MADM
must guarantee the monotonicity, i.e., the higher the score of each
attribute of the alternative is, the higher the ranking is. Meanwhile, decision-making
results obtained by this good method must be consistent with our intuitive judgment,
when dealing with the simplest problems, for which the decision making results can be obtained
by direct observation and comparison (see Examples~\ref{Exm-1} and \ref{Exm-2}).
In analyzing the TOPSIS method with crisp score values in $[0, 1]$, we can find that
the unit interval $[0, 1]$ has excellent algebraic and topological structures
as follows: (1) The unit interval $[0, 1]$ has a natural linear order structure $\leq$ and a
natural metric structure $| \cdot |$; (2) This natural metric structure $| \cdot |$
can ensure that the smaller the value in $[0, 1]$ is, the farther away from $1$,
and the closer away from $0$; (3) The order topology induced by the linear order
$\leq$ is consistent with the topology induced by the metric $| \cdot |$. These
good structures of $[0, 1]$ can ensure that the TOPSIS method of Hwang and Yoon~\cite{HY1981}
is monotonous. Recently, we~\cite{WWLCZ} proved that the space
$\tilde{\mathbb{I}}$ of all IFVs with the order topology induced by
the linear order $<_{_{\text{XY}}}$, defined in Definition~\ref{de-order(Xu)},
is not metrizable, i.e., there is no such good distance for $\tilde{\mathbb{I}}$
with the linear order $<_{_{\text{XY}}}$. Nevertheless, we still construct an
admissible distance with the order $\leq_{_{\text{XY}}}$ in \cite{WWLCZ}. In the
following, we will show that this distance is very important for our proposed
monotonous IF TOPSIS method.

\subsection{Limitation in TOPSIS method of Li~\cite{Li2014}}

First, we recall a fundamental IF TOPSIS method from~\cite{Li2014} and use two
examples show that the TOPSIS method in~\cite{Li2014} does not have
the basic monotonicity with the partial order $\subset$, which may lead the
decision-making result to be unreasonable and inconsistent with the actual situation,
when dealing with some even simplest decision-making problems.

The main process of IF TOPSIS method in~\cite[Section~3.3]{Li2014} is summarized as follows:

Step~1: Determine the alternatives $A=\{A_1, A_2, \ldots, A_n\}$
and attributes $\mathscr{O}=\{\mathscr{O}_1, \mathscr{O}_2, \ldots, \mathscr{O}_m\}$, respectively;

Step~2: Construct the IF decision matrix $R=(r_{ij}=\langle \mu_{ij}, \nu_{ij}\rangle)_{m\times n}$,
  as shown in Table~\ref{Tab-IFDM};

{Step~3: Determine the weights of the attributes expressed with the
IF weight vector $\omega=(\omega_1, \omega_2, \ldots,
  \omega_m)^{\top}$, where $\omega_j=\langle \rho_{j},
  \vartheta_j\rangle \in \tilde{\mathbb{I}}$;}

Step~4: Compute the weighted IF decision matrix $\overline{R}
=(\overline{r}_{ij}=\omega_{j}\otimes r_{ij})$ by Definition~\ref{Def-Int-Operations} (v),
i.e., $\overline{r}_{ij}=\langle \rho_{j} \mu_{ij}, \vartheta_j+\nu_{ij}-\vartheta_j\nu_{ij}\rangle$;

Step~5: Determine the IF positive ideal-point $\mathbf{A}^{+}=(\langle \mu^{+}_{1}, \nu^{+}_{1}\rangle,
  \langle \mu^{+}_{2}, \nu^{+}_{2}\rangle, \ldots, \langle \mu^{+}_{m}, \nu^{+}_{m}\rangle)^{\top}$
  and the IF negative ideal-point $\mathbf{A}^{-}=(\langle \mu^{-}_{1}, \nu^{-}_{1}\rangle,
  \langle \mu^{-}_{2}, \nu^{-}_{2}\rangle, \ldots, \langle \mu^{-}_{m}, \nu^{-}_{m}\rangle)^{\top}$ as follows:
  \begin{align*}
  & \mu^{+}_{j}=\max_{1\leq i\leq n}\{\bar{\mu}_{ij}\}, \quad
  \nu^{+}_{j}=\min_{1\leq i\leq n}\{\bar{\nu}_{ij}\}, \\
  & \mu^{-}_{j}=\min_{1\leq i\leq n}\{\bar{\mu}_{ij}\}, \quad
  \nu^{-}_{j}=\max_{1\leq i\leq n}\{\bar{\nu}_{ij}\};
 \end{align*}

Step~6: Compute the Euclidean distances $d_{\mathrm{Eu}}(A_i, \mathbf{A}^{+})$
  and $d_{\mathrm{Eu}}(A_i, \mathbf{A}^{-})$ of the alternatives
  $A_i$ ($i=1, 2, \ldots, n$) from $\mathbf{A}^{+}$ and $\mathbf{A}^{-}$
  by using formula~\eqref{Dis-Euclidean};

Step~7: Calculate the relative closeness degrees $\mathscr{C}_{i}$ of
  the alternatives $A_i$ ($i=1,2, \ldots,n$) to the IF positive
  ideal-point $\mathbf{A}^{+}$ by the following formula:
  $$
  \mathscr{C}_{i}=\frac{d_{\mathrm{Eu}}(A_i, \mathbf{A}^{-})}
  {d_{\mathrm{Eu}}(A_i, \mathbf{A}^{+})+d_{\mathrm{Eu}}(A_i, \mathbf{A}^{-})};
  $$

Step~8: Rank the alternatives
  $A_i$ ($i=1,2, \ldots,n$) according to the nonincreasing order
  of the relative closeness degrees $\mathscr{C}_{i}$ and
  select the most desirable alternative.

\begin{example}
\label{Exm-1}
Suppose that there exist $4$ alternatives $A_1$, $A_2$, $A_3$, $A_4$ evaluated with
respect to $2$ benefit attributes $\mathscr{O}_1$, $\mathscr{O}_2$. The sets of the
alternatives and attributes are denoted by $\{A_1, A_2, A_3, A_4\}$ and $\{\mathscr{O}_1,
\mathscr{O}_2\}$, respectively. Assume that the IF weight vector of $\mathscr{O}_1$
and $\mathscr{O}_2$ is $\omega=(\omega_1, \omega_2)^{\top}
=(\langle 1, 0\rangle, \langle 1, 0\rangle)^{\top}$. 
The IF decision-making matrix is expressed as shown in Table~\ref{Tab-2a}.
\begin{table}[H]	
	\centering
	\caption{IF decision matrix $R=(r_{ij})_{4\times 2}$}
	\label{Tab-2a}
     \scalebox{0.85}{
	\begin{tabular}{cccccc}
		\toprule
		$$ & $\mathscr{O}_{1}$ &  $\mathscr{O}_{2}$ \\
		\midrule
		$A_{1}$ & $\langle 0, 1\rangle$ &  $\langle 0, 1\rangle$ \\
		$A_{2}$ & $\langle 0.9, 0.01\rangle$ & $\langle 0.9, 0.01\rangle$ \\
        $A_{3}$ & $\langle 0.901, 0.007\rangle$ &  $\langle 0.901, 0.007\rangle$ \\
		$A_{4}$ & $\langle 1, 0\rangle$ & $\langle 1, 0\rangle$ \\
        \bottomrule
	\end{tabular}
      }
\end{table}
If we use the above TOPSIS method \cite[Section~3.3]{Li2014}, by direct calculation,
it can be verified that the weighted IF decision matrix is given as shown in Table~\ref{Tab-2b}.
\begin{table}[H]	
	\centering
	\caption{Weighted IF decision matrix $\overline{R}
=(\omega_{j}\otimes r_{ij})_{4\times 2}$}
	\label{Tab-2b}
     \scalebox{0.85}{
	\begin{tabular}{cccccc}
		\toprule
		$$ & $\mathscr{O}_{1}$ &  $\mathscr{O}_{2}$ \\
		\midrule
		$A_{1}$ & $\langle 0, 1\rangle$ &  $\langle 0, 1\rangle$ \\
		$A_{2}$ & $\langle 0.9, 0.01\rangle$ & $\langle 0.9, 0.01\rangle$ \\
        $A_{3}$ & $\langle 0.901, 0.007\rangle$ &  $\langle 0.901, 0.007\rangle$ \\
		$A_{4}$ & $\langle 1, 0\rangle$ & $\langle 1, 0\rangle$ \\
        \bottomrule
	\end{tabular}
      }
\end{table}
The IF positive ideal-point $\mathbf{A}^{+}$ and the IF
negative ideal-point $\mathbf{A}^{-}$ are obtained as follows:
$$
\mathbf{A}^{+}=(\langle 1, 0\rangle, \langle 1, 0\rangle) \text{ and }
\mathbf{A}^{-}=(\langle 0, 1 \rangle, \langle 0, 1 \rangle),
$$
respectively. According to the Euclidean distance of the alternatives $A_1$, $A_2$, $A_3$, and
$A_4$ from $\mathbf{A}^{+}$ and $\mathbf{A}^{-}$ obtained by~\cite[Eqs. (3.27) and (3.28)]{Li2014},
the relative closeness degrees $\mathscr{C}_{j}$ of the alternatives
$A_1$, $A_2$, $A_3$, and $A_4$ to the IF positive
ideal-point can be calculated as follows:
$$
\mathscr{C}_1=\frac{d_{\mathrm{Eu}}(A_1, \mathbf{A}^{-})}
{d_{\mathrm{Eu}}(A_1, \mathbf{A}^{+})+d_{\mathrm{Eu}}(A_1, \mathbf{A}^{-})}=0,
$$
$$
\mathscr{C}_2=\frac{d_{\mathrm{Eu}}(A_2, \mathbf{A}^{-})}
{d_{\mathrm{Eu}}(A_2, \mathbf{A}^{+})+d_{\mathrm{Eu}}(A_2, \mathbf{A}^{-})}=0.9085917,
$$
$$
\mathscr{C}_3=\frac{d_{\mathrm{Eu}}(A_3, \mathbf{A}^{-})}
{d_{\mathrm{Eu}}(A_3, \mathbf{A}^{+})+d_{\mathrm{Eu}}(A_3, \mathbf{A}^{-})}=0.9085194,
$$
and
$$
\mathscr{C}_4=\frac{d_{\mathrm{Eu}}(A_4, \mathbf{A}^{-})}
{d_{\mathrm{Eu}}(A_4, \mathbf{A}^{+})+d_{\mathrm{Eu}}(A_4, \mathbf{A}^{-})}=1,
$$
respectively. Therefore, the ranking order of
$A_1$, $A_2$, $A_3$, and $A_4$ is: $A_4\succ A_2 \succ A_3\succ A_1$. However,
$A_4\succ A_3 \succ A_2\succ A_1$ by a direct observation with $\langle 1, 0\rangle \supset
\langle 0.901, 0.007\rangle \supset \langle 0.9, 0.01\rangle \supset \langle 0, 1\rangle$.
This means that the ranking order obtained by the IF TOPSIS method
in \cite{Li2014} is not consistent with the real situation.
\end{example}

The following example demonstrates that the above IF TOPSIS may yield
some unreasonable decision-making results, even if we restrict the normalized weight vector
$\omega=(\omega_1, \omega_2, \ldots, \omega_m)^{\top}$ in
Step~3 to be positive real numbers, i.e., {$\omega_j\in (0, 1]$ and $\sum_{j=1}^{m}\omega_i=1$.}
\begin{example}
\label{Exm-2}
Suppose that there exist $4$ alternatives $A_1$, $A_2$, $A_3$, $A_4$ evaluated with respect to
$2$ benefit attributes $\mathscr{O}_1$, $\mathscr{O}_2$. The sets of the alternatives and attributes are
denoted by $\{A_1, A_2, A_3, A_4\}$ and $\{\mathscr{O}_1, \mathscr{O}_2\}$, respectively.
Assume that the weight vector of $\mathscr{O}_1$ and $\mathscr{O}_2$ is $\omega=(\omega_1, \omega_2)^{\top}
=(0.5, 0.5)^{\top}$.
 The IF decision-making matrix is expressed as shown in Table~\ref{Tab-1c}.
\begin{table}[H]	
	\centering
	\caption{IF decision-making matrix $R=(r_{ij})_{4\times 2}$}
	\label{Tab-1c}
     \scalebox{0.85}{
	\begin{tabular}{cccccc}
		\toprule
		$$ & $\mathscr{O}_{1}$ &  $\mathscr{O}_{2}$ \\
		\midrule
		$A_{1}$ & $\langle 0, 1\rangle$ &  $\langle 0, 1\rangle$ \\
		$A_{2}$ & $\langle 0.99, 0.0001\rangle$ & $\langle 0.99, 0.0001\rangle$ \\
        $A_{3}$ & $\langle 0.990199, 0.49\times 10^{-4}\rangle$ &  $\langle 0.990199, 0.49\times 10^{-4}\rangle$ \\
		$A_{4}$ & $\langle 1, 0\rangle$ & $\langle 1, 0\rangle$ \\
        \bottomrule
	\end{tabular}
      }
\end{table}
By direct calculation,
it can be verified that the weighted IF decision matrix is given
as shown in Table~\ref{Tab-2d}.
\begin{table}[H]	
	\centering
	\caption{Weighted IF decision matrix $\overline{R}
=(\omega_j\cdot r_{ij})_{4\times 2}$}
	\label{Tab-2d}
     \scalebox{0.85}{
	\begin{tabular}{cccccc}
		\toprule
		$$ & $\mathscr{O}_{1}$ &  $\mathscr{O}_{2}$ \\
		\midrule
		$A_{1}$ & $\langle 0, 1\rangle$ &  $\langle 0, 1\rangle$ \\
		$A_{2}$ & $\langle 0.9, 0.01\rangle$ & $\langle 0.9, 0.01\rangle$ \\
        $A_{3}$ & $\langle 0.901, 0.007\rangle$ &  $\langle 0.901, 0.007\rangle$ \\
		$A_{4}$ & $\langle 1, 0\rangle$ & $\langle 1, 0\rangle$ \\
        \bottomrule
	\end{tabular}
      }
\end{table}
If we use the TOPSIS method in \cite{Li2014}, by Example~\ref{Exm-1}, we know that the ranking order of
$A_1$, $A_2$, $A_3$, and $A_4$ is: $A_4\succ A_2 \succ A_3\succ A_1$. This is also an unreasonable
decision-making result.
\end{example}

\begin{remark}
(1) Examples~\ref{Exm-1} and \ref{Exm-2} show that, { for some weight vector
with either IFVs or real numbers, even for some simple decision-making problems, the TOPSIS method
in~\cite{Li2014} may lead to some unreasonable decision-making results.}

(2) Careful readers can verify that by applying the TOPSIS methods in
\cite{BGKA2009,BBM2012,MDJAR2019,RYU2020} to Examples~\ref{Exm-1}, the same result can
be obtained. This means that the TOPSIS methods in
\cite{BGKA2009,BBM2012,MDJAR2019,RYU2020} may produce unreasonable results when dealing
with the simplest decision-making problems.
\end{remark}

\subsection{Limitation in TOPSIS method of Chen et al.~\cite{CCL2016a}}

The above two examples show that the TOPSIS method in~\cite{Li2014} is not monotonous
with Atanassov's partial order $\subset$. Recently, Chen et al.~\cite{CCL2016a}
developed a monotonous TOPSIS method with the partial order $\subset$ based on
a new similarity measure. However, the following example shows that the
TOPSIS methods in \cite{CCL2016a} is not monotonous with the linear
order $\leq_{_\text{XY}}$.

The main process of IF TOPSIS method in~\cite{CCL2016a} is summarized as follows:

Step~1: Determine the alternatives $A=\{A_1, A_2, \ldots, A_n\}$
and attributes $\mathscr{O}=\{\mathscr{O}_1, \mathscr{O}_2, \ldots, \mathscr{O}_m\}$, respectively, and
construct the IF decision matrix $R=(r_{ij}=\langle \mu_{ij}, \nu_{ij}\rangle)_{m\times n}$,
as shown in Table~\ref{Tab-IFDM};

Step~2: Determine the IF positive ideal-point $\mathbf{A}^{+}=(\langle \mu^{+}_{1}, \nu^{+}_{1}\rangle,
  \langle \mu^{+}_{2}, \nu^{+}_{2}\rangle, \ldots, \langle \mu^{+}_{m}, \nu^{+}_{m}\rangle)^{\top}$
  and the IF negative ideal-point $\mathbf{A}^{-}=(\langle \mu^{-}_{1}, \nu^{-}_{1}\rangle,
  \langle \mu^{-}_{2}, \nu^{-}_{2}\rangle, \ldots, \langle \mu^{-}_{m}, \nu^{-}_{m}\rangle)^{\top}$ as follows:
  $$
\langle \mu^{+}_{j}, \nu^{+}_{j}\rangle=
\begin{cases}
\langle \max\limits_{1\leq i\leq n}\{{\mu}_{ij}\}, \min\limits_{1\leq i\leq n}\{{\nu}_{ij}\}\rangle, &
\mathscr{O}_{j}\in \mathscr{O}^{+}, \\
\langle \min\limits_{1\leq i\leq n}\{{\mu}_{ij}\}, \max\limits_{1\leq i\leq n}\{{\nu}_{ij}\}\rangle, &
\mathscr{O}_{j}\in \mathscr{O}^{-},
\end{cases}
$$
and
$$
\langle \mu^{-}_{j}, \nu^{-}_{j}\rangle=
\begin{cases}
\langle \min\limits_{1\leq i\leq n}\{{\mu}_{ij}\}, \max\limits_{1\leq i\leq n}\{{\nu}_{ij}\}\rangle, &
\mathscr{O}_{j}\in \mathscr{O}^{+}, \\
\langle \max\limits_{1\leq i\leq n}\{{\mu}_{ij}\}, \min\limits_{1\leq i\leq n}\{{\nu}_{ij}\}\rangle, &
\mathscr{O}_{j}\in \mathscr{O}^{-},
\end{cases}
$$
where $\mathscr{O}^{+}$ is the set of benefit attributes and
$\mathscr{O}^{-}$ is the set of cost attributes;

Step~3: Compute the degree of indeterminacy {$\pi_{j}^{+}=1-\mu^{+}_{j}-\nu^{+}_{j}$}
  of the positive ideal-point $\langle \mu^{+}_{j}, \nu^{+}_{j}\rangle$ for each attribute
  $\mathscr{O}_{j}$ ($j=1, 2, \ldots, n$);

Step~4: Compute the degree of indeterminacy {$\pi_{j}^{-}=1-\mu^{-}_{j}-\nu^{-}_{j}$}
  of the negative ideal-point $\langle \mu^{-}_{j}, \nu^{-}_{j}\rangle$ for each attribute
  $\mathscr{O}_{j}$ ($j=1, 2, \ldots, n$);

Step~5: Compute the degree of similarity $g_{ij}^{+}$ between the evaluating IFV
  $r_{ij}$ of the alternative $A_i$ with respect to the attribute $\mathscr{O}_{j}$ and the
  positive ideal-point $\langle \mu^{+}_{j}, \nu^{+}_{j}\rangle$ of
  the attribute $\mathscr{O}_{j}$ to construct the positive similarity matrix
  $G^{+}=(g_{ij}^{+})_{m\times n}$, where
  {$g_{ij}^{+}=1-\frac{|2(\mu_{j}^{+}-\mu_{ij})-(\nu_{j}^{+}-\nu_{ij})|}{3}\times
  (1-\frac{\pi_{j}^{+}+\pi_{ij}}{2})-\frac{|2(\nu_{j}^{+}-\nu_{ij})-
  (\mu_{j}^{+}-\mu_{ij})|}{3}\times \frac{\pi_{j}^{+}+\pi_{ij}}{2};$}

Step~6: Compute the degree of similarity $g_{ij}^{-}$ between the evaluating IFV
  $r_{ij}$ of the alternative $A_i$ with respect to the attribute $\mathscr{O}_{j}$ and the
  negative ideal-point $\langle \mu^{-}_{j}, \nu^{-}_{j}\rangle$ of
  the attribute $\mathscr{O}_{j}$ to construct the {negative} similarity matrix
  $G^{-}=(g_{ij}^{-})_{m\times n}$, where
  {$g_{ij}^{-}=1-\frac{|2(\mu_{j}^{-}-\mu_{ij})-(\nu_{j}^{-}-\nu_{ij})|}{3}\times
  (1-\frac{\pi_{j}^{-}+\pi_{ij}}{2})-\frac{|2(\nu_{j}^{-}-\nu_{ij})-
  (\mu_{j}^{-}-\mu_{ij})|}{3}\times \frac{\pi_{j}^{-}+\pi_{ij}}{2};$}

Step~7: Compute the weighted positive score $S_{i}^{+}=
  \sum_{j=1}^{m}\omega_{j} g_{ij}^{+}$ and
  the weighted negative score $S_{i}^{-}=\sum_{j=1}^{m}\omega_{j} g_{ij}^{-}$
  of each alternative $A_{i}$ ($i=1,2,\ldots, n$), where $\omega_j$ is the weight
  of criterion $\mathscr{O}_{j}$ such that $\omega_j\in (0, 1]$
  and $\sum_{j=1}^{m}\omega_j=1$;

Step~8: Compute the relative degree of closeness
   $T(A_i)=\frac{S_{i}^{+}}{S_{i}^{+}+S_{i}^{-}}$ of each alternative $A_{i}$.
   The larger the value of $T(A_i)$, the better the preference order of
   alternative $A_{i}$. Then, rank the alternatives
  $A_i$ ($i=1,2, \ldots,n$) according to the nonincreasing order
  of the relative closeness degrees $T(A_1)$, $T(A_2)$, $\ldots$,
  $T(A_{n})$.

\begin{example}
\label{Exm-3}
Suppose that there exist $4$ alternatives $A_1$, $A_2$, $A_3$, $A_4$ evaluated with respect to
$2$ benefit attributes $\mathscr{O}_1$, $\mathscr{O}_2$. The sets of the alternatives and attributes are
denoted by $\{A_1, A_2, A_3, A_4\}$ and $\{\mathscr{O}_1, \mathscr{O}_2\}$, respectively.
Assume that the weight vector of $\mathscr{O}_1$ and $\mathscr{O}_2$ is $\omega=(0.5, 0.5)^{\top}$. 
The IF decision-making matrix is expressed as shown in Table~\ref{Tab-2e}.
\begin{table}[H]	
	\centering
	\caption{IF decision matrix $R$}
	\label{Tab-2e}
     \scalebox{0.85}{
	\begin{tabular}{cccccc}
		\toprule
		$$ & $\mathscr{O}_{1}$ &  $\mathscr{O}_{2}$ \\
		\midrule
		$A_{1}$ & $\langle 0, 1\rangle$ &  $\langle 0, 1\rangle$ \\
		$A_{2}$ & $\langle 0.3, 0\rangle$ & $\langle 0.3, 0\rangle$ \\
        $A_{3}$ & $\langle 0.64, 0.36\rangle$ &  $\langle 0.64, 0.36\rangle$ \\
		$A_{4}$ & $\langle 1, 0\rangle$ & $\langle 1, 0\rangle$ \\
        \bottomrule
	\end{tabular}
      }
\end{table}
The IF positive ideal-point $\mathbf{A}^{+}$ and the IF
negative ideal-point $\mathbf{A}^{-}$ are obtained as follows:
$\mathbf{A}^{+}=(\langle 1, 0\rangle, \langle 1, 0\rangle)$ and
$\mathbf{A}^{-}=(\langle 0, 1 \rangle, \langle 0, 1 \rangle),$
respectively. By using  the TOPSIS method in~\cite{CCL2016a},
we obtain the positive similarity matrix $G^{+}$ and the negative similarity
matrix $G^{-}$ as follows:
$$
G^{+}=(g_{ij}^{+})_{4\times 2}=
\begin{bmatrix}
     0 & 0 \\

     0.615 & 0.615 \\

     0.64 & 0.64 \\

     1 & 1 \\
\end{bmatrix},
$$
and
$$
G^{-}=(g_{ij}^{-})_{4\times 2}=
\begin{bmatrix}
     1 & 1 \\

     0.385 & 0.385 \\

     0.36 & 0.36 \\

     0 & 0 \\
\end{bmatrix}.
$$
Then, the weighted positive scores $S^{+}_{i}=\omega_1g_{i1}^{+}+\omega_2g_{i2}^{+}$
($i=1,2,3,4$) and the weighted negative scores
${S_{i}^{-}}=\omega_1g_{i1}^{-}+\omega_2g_{i2}^{-}$ ($i=1,2,3,4$)
of the alternatives $A_1$, $A_2$, $A_3$, and $A_4$ can be calculated as follows:
$$
S_{1}^{+}=0, \ S_{2}^{+}=0.615, \ S_{3}^{+}=0.64, \ S_{4}^{+}=1,
$$
and
$$
S_{1}^{-}=1, \ S_{2}^{+}=0.385, \ S_{3}^{+}=0.36, \ S_{4}^{-}=0.
$$
Therefore, the relative degree of closeness $T(A_{i})=\frac{S_{i}^{+}}{S_{i}^{+}+S_{i}^{-}}$
($i=1,2,3,4$) of the alternatives
$A_1$, $A_2$, $A_3$, and $A_4$ are given as follows:
$$
T(A_{1})=0, \ T(A_{2})=0.615, \  T(A_{3})=0.64, \  T(A_{4})=1,
$$
and thus the ranking order of $A_1$, $A_2$, $A_3$, and $A_4$ is:
$A_4\succ A_3 \succ A_2\succ A_1$. However, it can be verified that
$A_4\succ A_2 \succ A_3 \succ A_1$
by a direct observation with $\langle 1, 0\rangle \geq_{_{\text{XY}}}
\langle 0.3, 0\rangle \geq_{_{\text{XY}}} \langle 0.64, 0.36\rangle \geq_{_{\text{XY}}} \langle 0, 1\rangle$.
\end{example}

Summing up Examples~\ref{Exm-1}--\ref{Exm-3}, an interesting question is {\it whether there exists
an IF TOPSIS method that is monotonous with the linear order $\leq_{_\text{XY}}$ or $\leq_{_\text{ZX}}$}?
In the following section, we will establish an IF TOPSIS method that is monotonous with the
linear order $\leq_{_\text{XY}}$ or $\leq_{_\text{ZX}}$.

\section{A monotonous IF TOPSIS method}\label{Sec-6}

The fundamental cause for counterintuitive decision-making results in
Examples~\ref{Exm-1}--\ref{Exm-3} lies in the structure of metrics for IFVs.
In \cite{WWLCZ}, we defined a metric $\varrho$ in $\tilde{\mathbb{I}}$ as follows:
for $\alpha$, $\beta\in \tilde{\mathbb{I}}$,
\begin{equation*}
\varrho(\alpha, \beta)=
\begin{cases}
\frac{1}{3}(1+|s(\alpha)-s(\beta)|), & s(\alpha)\neq s(\beta), \\
\frac{1}{3}(|h(\alpha)-h(\beta)|), & s(\alpha)=s(\beta),
\end{cases}
\end{equation*}
where $s(\alpha)$ and $h(\alpha)$ are the score degree and the accuracy degree of $\alpha$, respectively.
Furthermore, we~\cite{WWLCZ} proved the following basic properties of $\varrho$.
\begin{theorem}[{\textrm{\protect\cite{WWLCZ}}}]
\label{Thm-metric}
\begin{enumerate}[{\rm (1)}]
  \item $\varrho(\alpha, \beta)\in [0, 1]$ and $\varrho(\alpha, \beta)=0$ if and only if $\alpha=\beta$.
  \item $\varrho(\alpha, \beta)=1$ if and only if ($\alpha=\langle 0, 1\rangle$ and
$\beta=\langle 1, 0\rangle$) or ($\alpha=\langle 1, 0\rangle$ and
$\beta=\langle 0, 1\rangle$).
  \item $\varrho(\alpha, \beta)=\varrho(\beta, \alpha)$.
  \item For any $\alpha$, $\beta$, $\gamma \in \tilde{\mathbb{I}}$,
$\varrho(\alpha, \beta)+\varrho(\beta, \gamma)\geq \varrho(\alpha, \gamma)$.
  \item For any $\alpha$, $\beta$, $\gamma \in \tilde{\mathbb{I}}$, if $\alpha\leq_{_{\text{XY}}}
\beta \leq_{_{\text{XY}}} \gamma$, then $\varrho(\alpha, \beta)\leq \varrho(\alpha, \gamma)$
and $\varrho(\beta, \gamma)\leq \varrho(\alpha, \gamma)$.
\end{enumerate}
\end{theorem}

{Based on the similarity function $L(\alpha)$, similarly to the metric $\varrho$,
define the parametric metrics $\varrho^{(\lambda)}$ and $\tilde{\varrho}^{(\lambda)}$
in $\tilde{\mathbb{I}}$ as follows: for $\alpha$, $\beta\in \tilde{\mathbb{I}}$,
\begin{equation}
\label{rho-equ}
\begin{split}
& \varrho^{(\lambda)}(\alpha, \beta)\\
= &
\begin{cases}
\frac{1}{1+2\lambda}(1+\lambda \cdot |s(\alpha)-s(\beta)|), & s(\alpha)\neq s(\beta), \\
\frac{1}{1+2\lambda}(|h(\alpha)-h(\beta)|), & s(\alpha)=s(\beta),
\end{cases}
\end{split}
\end{equation}
and
\begin{equation}
\label{rho-equ-1}
\begin{split}
& \tilde{\varrho}^{(\lambda)}(\alpha, \beta)\\
=&
\begin{cases}
\frac{1}{1+\lambda}(1+\lambda\cdot |L(\alpha)-L(\beta)|), & L(\alpha)\neq L(\beta), \\
\frac{1}{1+\lambda}(|h(\alpha)-h(\beta)|), & L(\alpha)=L(\beta),
\end{cases}
\end{split}
\end{equation}
where $\lambda\geq 1$ is a parameter, and $L(\alpha)$ and $h(\alpha)$ are the similarity function and the accuracy
degree of $\alpha$, respectively.

Let $A$ and $B$ be two aggregation functions satisfying the condition in
Proposition~\ref{Prop-Linear-Order}. Based on Proposition~\ref{Prop-Linear-Order},
define another parametric metric $\varrho_{_{A, B}}^{(\lambda)}$ in $\tilde{\mathbb{I}}$ as follows:
for $\alpha$, $\beta\in \tilde{\mathbb{I}}$,
\begin{equation}
\label{rho-equ-2}
\begin{split}
&\varrho_{_{A, B}}^{(\lambda)}(\alpha, \beta)\\
=&
\begin{cases}
\frac{1}{2}(1+|\overline{A}(\alpha)-\overline{A}(\beta)|), & \overline{A}(\alpha)\neq \overline{A}(\beta), \\
\frac{1}{2}(|\overline{B}(\alpha)-\overline{B}(\beta)|), & \overline{A}(\alpha)=\overline{A}(\beta),
\end{cases}
\end{split}
\end{equation}
where $\lambda\geq 1$ is a parameter, $\overline{A}(\alpha)=A(\mu_{\alpha}, 1-\nu_{\alpha})$, and
$\overline{B}(\alpha)=B(\mu_{\alpha}, 1-\nu_{\alpha})$. In particular,
by taking $A=K_{\gamma_1}$ and $B=K_{\gamma_2}$ with $\gamma_1\neq \gamma_2$
and direct calculation, we have
\begin{equation}
\label{rho-equ-3}
\begin{split}
&\varrho_{_{K_{\gamma_1}, K_{\gamma_2}}}^{(\lambda)}(\alpha, \beta)\\
=&
\begin{cases}
\frac{1}{1+\lambda}(1+\lambda\cdot E_{\gamma_1}(\alpha, \beta)), & E_{\gamma_1}(\alpha, \beta)\neq 0, \\
\frac{1}{1+\lambda}E_{\gamma_2}(\alpha, \beta), & E_{\gamma_1}(\alpha, \beta)=0,
\end{cases}
\end{split}
\end{equation}
where $E_{\gamma}(\alpha, \beta)=|\overline{K}_{\gamma}(\alpha)-
\overline{K}_{\gamma}(\beta)|=|(1-\gamma)(\mu_{\alpha}-\mu_{\beta})
-\gamma(\nu_{\alpha}-\nu_{\beta})|$.

\begin{remark}
The parameter $\gamma_1$ in Eq.~\eqref{rho-equ-3} can be regarded as the preference for decision-makers
to choose the membership and non-membership:

(1) If $\gamma_1>0.5$, then the decision-makers prefer non-membership to membership,
i.e., the decision-makers are pessimistic.

(2) If $\gamma_1<0.5$, then the decision-makers prefer membership to non-membership,
i.e., the decision-makers are optimistic.

(3) If $\gamma_1=0.5$, then the decision-makers have no preference for membership
and non-membership, i.e., the decision-makers are neutral.
\end{remark}

Similarly to the proof of Theorem~\ref{Thm-metric} in \cite{WWLCZ},
we can prove that the metrics $\varrho^{(\lambda)}$, $\tilde{\varrho}^{(\lambda)}$,
and $\varrho_{_{A,B}}^{(\lambda)}$ have the following basic properties.
\begin{theorem}
\label{Thm-metric-1}
Let $\lambda\geq 1$ and $\rho\in \{\varrho^{(\lambda)}, \tilde{\varrho}^{(\lambda)},
\varrho_{_{A,B}}^{(\lambda)}\}$. Then,
\begin{enumerate}[{\rm (1)}]
  \item $\rho(\alpha, \beta)\in [0, 1]$ and $\rho(\alpha, \beta)=0$ if and only if $\alpha=\beta$.
  \item $\rho(\alpha, \beta)=1$ if and only if ($\alpha=\langle 0, 1\rangle$ and
$\beta=\langle 1, 0\rangle$) or ($\alpha=\langle 1, 0\rangle$ and
$\beta=\langle 0, 1\rangle$).
  \item $\rho(\alpha, \beta)=\varrho(\beta, \alpha)$.
  \item For any $\alpha$, $\beta$, $\gamma \in \tilde{\mathbb{I}}$,
$\rho(\alpha, \beta)+\rho(\beta, \gamma)\geq \rho(\alpha, \gamma)$.
  \item For any $\alpha$, $\beta$, $\gamma \in \tilde{\mathbb{I}}$, if $\alpha \leq
\beta \leq \gamma$, then $\rho(\alpha, \beta)\leq \rho(\alpha, \gamma)$
and $\rho(\beta, \gamma)\leq \rho(\alpha, \gamma)$.
\end{enumerate}
\end{theorem}}

For the MADM problem with IFSs,
by using the three metrics defined by Eqs.~\eqref{rho-equ}--\eqref{rho-equ-2}, we propose a new
IF TOPSIS method as follows:

Step~1:  (Construct the decision matrix)
  Supposing that the decision-maker gave the rating (or evaluation) of each
  alternative $A_i\in A$ ($i=1, 2, \ldots, n$) on each attribute $\mathscr{O}_j$
  ($j=1, 2, \ldots, m$) in the form of IFNs $r_{ij}=\langle \mu_{ij},
  \nu_{ij}\rangle$, construct an IF decision matrix $R=(r_{ij})_{n\times m}$
  as shown in Table~\ref{Tab-IFDM}.

Step~2: (Normalize the decision matrix) Transform the IF decision
matrix $R=(r_{ij})_{n\times m}$ to the normalized IF decision
matrix $\overline{R}=(\bar{r}_{ij})_{n\times m}=(\langle \bar{\mu}_{ij},
\bar{\nu}_{ij}\rangle)_{n\times m}$ as follows:
$$
\bar{r}_{ij}=
\begin{cases}
r_{ij}, & \text{for benefit attribute } \mathscr{O}_{j}, \\
r_{ij}^{\complement}, & \text{for cost attribute } \mathscr{O}_{j},
\end{cases}
$$
where $r_{ij}^{\complement}$ is the complement of {$r_{ij}$}.

Step~3: (Determine the positive and negative ideal-points)
  Determine the
  IF positive ideal-point $\mathbf{A}^{+}=(\langle \mu^{+}_{1}, \nu^{+}_{1}\rangle,
  \langle \mu^{+}_{2}, \nu^{+}_{2}\rangle, \ldots, \langle \mu^{+}_{m}, \nu^{+}_{m}\rangle)^{\top}$
  and IF negative ideal-point $\mathbf{A}^{-}=(\langle \mu^{-}_{1}, \nu^{-}_{1}\rangle,
  \langle \mu^{-}_{2}, \nu^{-}_{2}\rangle, \ldots, \langle \mu^{-}_{m}, \nu^{-}_{m}\rangle)^{\top}$ as follows:
  \begin{align*}
  & \mu^{+}_{j}=\max_{1\leq i\leq n}\{\bar{\mu}_{ij}\}, \quad
  \nu^{+}_{j}=\min_{1\leq i\leq n}\{\bar{\nu}_{ij}\},\\
  & \mu^{-}_{j}=\min_{1\leq i\leq n}\{\bar{\mu}_{ij}\}, \quad
  \nu^{-}_{j}=\max_{1\leq i\leq n}\{\bar{\nu}_{ij}\}.
  \end{align*}

Step~4: (Compute the weighted similarity measures)
  Choose $\lambda\geq 1$ and compute the weighted similarity measures between
  the alternatives $A_i$ ($i=1,2, \ldots,n$) and the IF positive
  ideal-point $\mathbf{A}^+$ and between
  the alternatives $A_i$ ($i=1,2, \ldots,n$) and the IF negative
  ideal-point $\mathbf{A}^{-}$ by using the following formulas:
  {\small \begin{align}
  & \mathbf{S}(A_i, \mathbf{A}^{+})=1-\sum_{j=1}^{m}\omega_{j}\cdot \varrho^{(\lambda)}(\langle \bar{\mu}_{ij},
\bar{\nu}_{ij}\rangle, \langle \mu^{+}_{j}, \nu^{+}_{j}\rangle),  \label{Eq-1} \\
 (\text{resp., } & \mathbf{S}(A_i, \mathbf{A}^{+})=1-\sum_{j=1}^{m}\omega_{j}\cdot \tilde{\varrho}^{(\lambda)}(\langle \bar{\mu}_{ij},
\bar{\nu}_{ij}\rangle, \langle \mu^{+}_{j}, \nu^{+}_{j}\rangle),  \label{Eq-1-a} \\
& \mathbf{S}(A_i, \mathbf{A}^{+})=1-\sum_{j=1}^{m}\omega_{j}\cdot \varrho_{_{A, B}}^{(\lambda)}(\langle \bar{\mu}_{ij},
\bar{\nu}_{ij}\rangle, \langle \mu^{+}_{j}, \nu^{+}_{j}\rangle)),  \label{Eq-1-aa}
  \end{align}
  and
  \begin{align}
  & \mathbf{S}(A_i, \mathbf{A}^{-})=1-\sum_{j=1}^{m}\omega_{j}\cdot \varrho^{(\lambda)}(\langle \bar{\mu}_{ij},
\bar{\nu}_{ij}\rangle, \langle \mu^{-}_{j}, \nu^{-}_{j}\rangle),  \label{Eq-2}\\
  (\text{resp., } & \mathbf{S}(A_i, \mathbf{A}^{-})=1-\sum_{j=1}^{m}\omega_{j}\cdot \tilde{\varrho}^{(\lambda)}(\langle \bar{\mu}_{ij},
\bar{\nu}_{ij}\rangle, \langle \mu^{-}_{j}, \nu^{-}_{j}\rangle),   \label{Eq-2-b}\\
& \mathbf{S}(A_i, \mathbf{A}^{-})=1-\sum_{j=1}^{m}\omega_{j}\cdot \varrho_{_{A, B}}^{(\lambda)}(\langle \bar{\mu}_{ij},
\bar{\nu}_{ij}\rangle, \langle \mu^{-}_{j}, \nu^{-}_{j}\rangle)).   \label{Eq-2-bb}
\end{align}}
By Theorems~\ref{Thm-metric}
and \ref{Thm-metric-1}, it is easy to see that the similarity measures
obtained by Eqs.~\eqref{Eq-1}--\eqref{Eq-1-aa} are admissible similarity measures
with the orders $\leq_{_{\text{XY}}}$, $\leq_{_{\text{ZX}}}$, and
$\leq_{_{A,B}}$, respectively.

  Step~5: (Compute the relative closeness degrees)
  Calculate the relative closeness degrees $\mathscr{C}_{i}$ of
  the alternatives $A_i$ ($i=1,2, \ldots,n$) to the IF positive
  ideal-point $\mathbf{A}^+$ by using the following formula:
  \begin{equation}
  \label{Eq-3}
  \mathscr{C}_{i}=\frac{\mathbf{S}(A_i, \mathbf{A}^{+})}
  {\mathbf{S}(A_i, \mathbf{A}^{+})+\mathbf{S}(A_i, \mathbf{A}^{-})}.
  \end{equation}

  Step~6: (Rank the alternatives) Rank the alternatives
  $A_i$ ($i=1,2, \ldots,n$) according to the nonincreasing order
  of the relative closeness degrees $\mathscr{C}_{i}$ and
  select the most desirable alternative.

\begin{remark}
(1) By Theorems~\ref{Thm-metric} and \ref{Thm-metric-1}, it is easy to see that
$\mathbf{S}(A_i, \mathbf{A}^{+})+\mathbf{S}(A_i, \mathbf{A}^{-})$ in Eq.~\eqref{Eq-3}
is always nonzero. This overcomes the limitation that many TOPSIS method
may lead to the situation that the denominator is equal to $0$ when computing
the relative closeness degrees.

{(2) For practical MADM problems, in order to eliminate the effect of the constant term
$1$ in formulas~\eqref{rho-equ}--\eqref{rho-equ-3} as much as possible, the parameter
$\lambda$ should be chosen as large as possible.}
\end{remark}

\begin{theorem}[\textrm{Monotonicity}]
\label{Mono-Thm}
Using Eqs.~\eqref{Eq-1} and \eqref{Eq-2},
the above proposed method is increasing with the linear order $\leq_{_\text{XY}}$, i.e.,
for the MADM problem expressed in Table~\ref{Tab-IFDM}, if there exist $1\leq i_1, i_2 \leq n$
such that $\bar{r}_{i_1j}\leq _{_{\text{XY}}}\bar{r}_{i_2j}$ holds for all $1\leq j\leq m$,
then $\mathscr{C}_{i_1}\leq \mathscr{C}_{i_2}$, i.e., $A_{i_2}$ is better than $A_{i_1}$
ranked by the proposed method. In particular, the proposed method is increasing with
Atanassov's order `$\subset$'.
\end{theorem}

\begin{IEEEproof}
Fix $\lambda\geq 1$. Let $\mathbf{A}^{+}=(\langle \mu^{+}_{1}, \nu^{+}_{1}\rangle,
  \langle \mu^{+}_{2}, \nu^{+}_{2}\rangle, \ldots, $ $\langle \mu^{+}_{m}, \nu^{+}_{m}\rangle)^{\top}$
  and  $\mathbf{A}^{-}=(\langle \mu^{-}_{1}, \nu^{-}_{1}\rangle,
  \langle \mu^{-}_{2}, \nu^{-}_{2}\rangle, \ldots, \langle \mu^{-}_{m}, \nu^{-}_{m}\rangle)^{\top}$
  be the IF positive ideal-point and the IF negative ideal-point obtained by Step~3, respectively.
  Clearly, $\langle \mu^{-}_{j}, \nu^{-}_{j} \rangle\leq_{_{\text{XY}}}\bar{r}_{i_1j}\leq _{_{\text{XY}}}
\bar{r}_{i_2j}\leq_{_{\text{XY}}} \langle \mu^{+}_{j}, \nu^{+}_{j}\rangle$. By Theorem~\ref{Thm-metric-1},
we have
$$
\varrho^{(\lambda)}(\langle \mu^{-}_{j}, \nu^{-}_{j} \rangle, \bar{r}_{i_1j})\leq
\varrho^{(\lambda)}(\langle \mu^{-}_{j}, \nu^{-}_{j} \rangle, \bar{r}_{i_2j}),
$$
and
$$
\varrho^{(\lambda)}(\langle \mu^{+}_{j}, \nu^{+}_{j}\rangle, \bar{r}_{i_2j})
\leq \varrho^{(\lambda)}(\langle \mu^{+}_{j}, \nu^{+}_{j}\rangle, \bar{r}_{i_1j}),
$$
and thus,
$$
\mathbf{S}(A_{i_1}, \mathbf{A}^{-})\geq \mathbf{S}(A_{i_2}, \mathbf{A}^{-}),
$$
and
$$
\mathbf{S}(A_{i_1}, \mathbf{A}^{+})\leq \mathbf{S}(A_{i_2}, \mathbf{A}^{+}) \text{ by Eqs.~\eqref{Eq-1}
 and \eqref{Eq-2}}.
$$
This, together with Eq.~\eqref{Eq-3}, implies that

(1) if $\mathbf{S}(A_{i_1}, \mathbf{A}^{+})=0$, then
$$
\mathscr{C}_{i_1}=\frac{\mathbf{S}(A_{i_1}, \mathbf{A}^{+})}
{\mathbf{S}(A_{i_1}, \mathbf{A}^{+})+\mathbf{S}(A_{i_1},\mathbf{A}^{-})}
=0\leq \mathscr{C}_{i_2};
$$

(2) if $\mathbf{S}(A_{i_1}, \mathbf{A}^{+})>0$, then $\mathbf{S}(A_{i_2}, \mathbf{A}^{+})\geq
\mathbf{S}(A_{i_1}, \mathbf{A}^{+})>0$, and thus
\begin{align*}
\mathscr{C}_{i_1}=& \frac{\mathbf{S}(A_{i_1}, \mathbf{A}^{+})}{\mathbf{S}(A_{i_1}, \mathbf{A}^{+})
+\mathbf{S}(A_{i_1}, \mathbf{A}^{-})}
=\frac{1}{1+\frac{\mathbf{S}(A_{i_1}, \mathbf{A}^{-})}{\mathbf{S}(A_{i_1}, \mathbf{A}^{+})}}\\
\leq &
\frac{1}{1+\frac{\mathbf{S}(A_{i_2}, \mathbf{A}^{-})}{\mathbf{S}(A_{i_2}, \mathbf{A}^{+})}}=
\frac{\mathbf{S}(A_{i_2}, \mathbf{A}^{+})}{\mathbf{S}(A_{i_2}, \mathbf{A}^{+})+\mathbf{S}(A_{i_2},
\mathbf{A}^{-})}=
\mathscr{C}_{i_2}.
\end{align*}
Therefore, $\mathscr{C}_{i_1}\leq \mathscr{C}_{i_2}$.
\end{IEEEproof}

By Theorem~\ref{Thm-metric-1}, similarly to the proof of
Theorem~\ref{Mono-Thm}, it is not difficult to check that
the following result holds.

{\begin{theorem}[\textrm{Monotonicity}]
Using Eqs.~\eqref{Eq-1-a} and \eqref{Eq-2-b} (resp.,
Eqs.~\eqref{Eq-1-aa} and \eqref{Eq-2-bb}), the above proposed method is increasing
with the linear order $\leq_{_\text{ZX}}$ (resp., $\leq_{_{A, B}}$).
\end{theorem}}

{\begin{example}[\textrm{Continuation of Example~\ref{Exm-2}}]
\label{Exm-1-Cont}
Consider the MADM problem described in Example~\ref{Exm-2}. If the proposed
TOPSIS method in this section is used based on Eq.~\eqref{rho-equ} with $\lambda=1$,
by direct calculation, it can be verified that
$\mathscr{C}_1=0$, $\mathscr{C}_2=0.99495$, $\mathscr{C}_3=0.995075$, $\mathscr{C}_4=1$.
Thus, the ranking order of the alternatives $A_1$, $A_2$, $A_3$, and $A_4$ is:
$A_4\succ A_3\succ A_2\succ A_1$, which is consistent with the result
obtained by directly observing in Example~\ref{Exm-2}.
\end{example}}

{\begin{example}[\textrm{Continuation of Example~\ref{Exm-3}}]
\label{Exm-2-Cont}
Consider the MADM problem described in Example~\ref{Exm-3}. If the proposed
TOPSIS method in this section is used based on Eq.~\eqref{rho-equ} with $\lambda=1$,
by direct calculation, it can be verified that
$\mathscr{C}_1=0$, $\mathscr{C}_2=0.65$, $\mathscr{C}_3=0.64$, $\mathscr{C}_4=1$.
Thus, the ranking order of the alternatives $A_1$, $A_2$, $A_3$, and $A_4$ is:
$A_4\succ A_2\succ A_3\succ A_1$, which is consistent with the result
obtained by directly observing in Example~\ref{Exm-3}.
\end{example}

\begin{remark}
Observing from Examples~\ref{Exm-1-Cont} and \ref{Exm-2-Cont}, it can be seen that
the proposed TOPSIS method can effectively overcome the limitations of the TOPSIS
methods in \cite{Li2014,CCL2016a}, which is consistent with the result proved
in Theorem~\ref{Mono-Thm}. Furthermore, this shows that the proposed TOPSIS
method is superior to those in \cite{Li2014,CCL2016a}.
\end{remark}}

\section{Illustrative examples}\label{Sec-7}

This section provides two practical examples to illustrate the efficiency
of the above proposed TOPSIS method. One is an IF MADM problem on the choice
of suppliers in the supply chain management (see Example~\ref{Exm-WuX}).
The ranking order obtained by the proposed TOPSIS method is slightly
different from the results obtained by those TOPSIS methods in~\cite{AAB2021,
BG2016,ZZY2020,SMLXC2018}. However, the most desirable alternatives are
consistent. The other is an IF MADM problem on the choice of project managers
(see Example~\ref{Exm-WuX-1}). The ranking order, obtained by the proposed
TOPSIS method under the case that the decision-maker is neutral or pessimistic,
is consistent with those results obtained by the TOPSIS methods
in~\cite{CCL2016a,WW2008,AAB2021,BG2016,ZZY2020}.

\begin{example}
[{\textrm{\protect\cite[Example 5.1]{CCL2016a}}}]
\label{Exm-WuX}
Assume that there are five alternatives $A_1$, $A_2$, $A_3$, $A_4$, and
$A_5$ of suppliers and four attributes $\mathscr{O}_1$, $\mathscr{O}_2$,
$\mathscr{O}_3$, and $\mathscr{O}_4$ to assess these five alternatives,
so as to choose the best supplier among these five alternatives in the supply chain management,
where $\mathscr{O}_1$ is the ``Product Quality", $\mathscr{O}_2$ is the ``Service",
$\mathscr{O}_3$ is the ``Delivery", $\mathscr{O}_4$ is the ``Sustainability" and $\mathscr{O}_1$,
$\mathscr{O}_2$, $\mathscr{O}_3$, and $\mathscr{O}_4$ are
benefit attributes, with weight vector $\omega=(0.25, 0.4, 0.2, 0.15)^{\top}$.

Step~1:  (Construct the decision matrix)
The decision matrix $R=(r_{ij})_{5\times 4}$ given by the decision maker is listed
in Table~\ref{Tab-3}.
\begin{table}[H]
\caption{The decision matrix $R$}
\label{Tab-3}\centering
\scalebox{0.85}{
	\begin{tabular}{cccccc}
		\toprule
		$$ & $\mathscr{O}_{1}$ &  $\mathscr{O}_{2}$ & $\mathscr{O}_{3}$ & $\mathscr{O}_{4}$ \\
		\midrule
		$A_{1}$ & $\langle0.6, 0.3\rangle$ & $\langle0.5, 0.2\rangle$ &
           $\langle0.2, 0.5\rangle$ & $\langle 0.1, 0.6\rangle$\\
		$A_{2}$ & $\langle0.8, 0.2\rangle$ & $\langle0.8, 0.1\rangle$ &
           $\langle0.6, 0.1\rangle$ & $\langle0.3, 0.4\rangle$ \\
        $A_{3}$ & $\langle0.6, 0.3\rangle$ & $\langle0.4, 0.3\rangle$ &
           $\langle0.4, 0.2\rangle$ & $\langle0.5, 0.2\rangle$ \\
        $A_{4}$ & $\langle0.9, 0.1\rangle$ & $\langle0.5, 0.2\rangle$ &
           $\langle0.2, 0.3\rangle$ & $\langle0.1, 0.5\rangle$ \\
        $A_{5}$ & $\langle0.7, 0.1\rangle$ & $\langle0.3, 0.2\rangle$ &
           $\langle0.6, 0.2\rangle$ & $\langle0.4, 0.2\rangle$ \\
        \bottomrule
	\end{tabular}
      }
\end{table}

  Step~2: (Normalize the decision matrix) Since $\mathscr{O}_1$,
  $\mathscr{O}_2$, $\mathscr{O}_3$, and $\mathscr{O}_4$ are
all benefit attributes, we have $\overline{R}=(\bar{r}_{ij})_{5\times 4}=R$.

  Step~3: (Determine the positive and negative ideal-points)
  The IF positive ideal-point is
  $$
  \mathbf{A}^{+}=(\langle 0.9 , 0.1\rangle,
  \langle 0.8, 0.1\rangle, \langle 0.6, 0.1\rangle, \langle 0.5, 0.2\rangle)^{\top},
  $$
  and IF negative ideal-point is
  $$
  \mathbf{A}^{-}=(\langle 0.6, 0.3\rangle,
  \langle 0.3, 0.3\rangle, \langle 0.2, 0.5\rangle, \langle 0.1, 0.6\rangle)^{\top}.
  $$

  Steps~4 and 5: (Compute the relative closeness degrees)
  Choose $\lambda=100$ and calculate the relative closeness degrees $\mathscr{C}_{i}$ of
  the alternatives $A_i$ ($i=1, 2, 3, 4, 5$) to the IF positive
  ideal-point $\mathbf{A}^+$ by Eqs.~\eqref{Eq-1}, \eqref{Eq-2}, and
  \eqref{Eq-3}: $\mathscr{C}_{1}=0.4321$, $\mathscr{C}_{2}=0.5709$,
  $\mathscr{C}_{3}=0.4750$, $\mathscr{C}_{4}=0.4876$,
  $\mathscr{C}_{5}=0.5053.$

  Step~6: (Rank the alternative) Because
  $\mathscr{C}_2> \mathscr{C}_5> \mathscr{C}_4> \mathscr{C}_3
  >\mathscr{C}_1$, the ranking order of the alternatives $A_i$ ($i=1, 2, 3, 4, 5$) is:
   $A_2\succ A_5\succ A_4\succ A_3
  \succ A_1.$

Repeating Steps~1--3, by applying Eqs.~\eqref{Eq-1-a} and \eqref{Eq-2-b}, we obtain the following result:

Step~4 and 5: (Compute the relative closeness degrees)
  Choose $\lambda=100$ and calculate the relative closeness degrees $\mathscr{C}_{i}$ of
  the alternatives $A_i$ ($i=1, 2, 3, 4, 5$) to the IF positive
  ideal-point $\mathbf{A}^+$ by Eqs.~\eqref{Eq-1-a}, \eqref{Eq-2-b}, and
  \eqref{Eq-3}:
  $\mathscr{C}_{1}=0.4371$, $\mathscr{C}_{2}=0.5618$, $\mathscr{C}_{3}=0.4694$,
  $\mathscr{C}_{4}=0.4922$, $\mathscr{C}_{5}=0.4925.$

  Step~6: (Rank the alternatives) Because
  $\mathscr{C}_2> \mathscr{C}_5> \mathscr{C}_4> \mathscr{C}_3
  > \mathscr{C}_1$, the ranking order of the alternatives $A_i$
  ($i=1, 2, 3, 4, 5$) is: $A_2\succ A_5\succ A_4\succ A_3
  \succ A_1.$
\end{example}

\subsection*{Comparative analysis}

{From Table~\ref{Tab-Exm-WuX}, which shows a comparison of the ranking orders of the alternatives
in Example~\ref{Exm-WuX} for different MADM methods, it can be observed that (1) our results are exactly
the same, which are consistent with the result $A_2\succ A_5\succ A_4\succ A_3 \succ A_1$
in \cite{CCL2016a,WW2008,ZCK2019}; (2) the result obtained by Xu's IFWA operator in
\cite{Xu2007} is different from the results obtained by all other methods;
(3) since the TOPSIS method of B{\"u}y{\"u}k{\"o}zkan and G{\"u}lery{\"u}z in \cite{BG2016}
is based on the normalized Euclidean distance defined by Eq.~\eqref{Dis-Euclidean}, which
does not satisfy the axiomatic definition of IF distance measure
(see Example~\ref{Exm-No-Metric-1}), the ranking result
$A_2\succ A_5\succ A_3\succ A_4\succ A_1$ may be unreasonable;
(4) the best choice is always $A_2$.
\begin{table}[H]	
	\centering
	\caption{A comparison of the ranking orders of the alternatives
in Example~\ref{Exm-WuX} for different MADM methods}
	\label{Tab-Exm-WuX}
     \scalebox{0.85}{
	\begin{tabular}{cccccc}
		\toprule
		Methods &  Ranking orders \\
		\midrule
		Chen et al.'s TOPSIS method in \cite{CCL2016a} & $A_2\succ A_5\succ A_4\succ A_3\succ A_1$ \\
		Wang and Wei's TOPSIS method in \cite{WW2008} & $A_2\succ A_5\succ A_4\succ A_3\succ A_1$ \\
        Altan~Koyuncu et al.'s TOPSIS method in \cite{AAB2021} & $A_2\succ A_5\succ A_3\succ A_4\succ A_1$ \\
        B{\"u}y{\"u}k{\"o}zkan and G{\"u}lery{\"u}z's TOPSIS method in \cite{BG2016} & $A_2\succ A_5\succ A_3\succ A_4\succ A_1$ \\
        Zhang et al.'s TOPSIS method in \cite{ZZY2020} & $A_2\succ A_5\succ A_3\succ A_4\succ A_1$ \\
        Zeng et al.'s VIKOR method in \cite{ZCK2019} & $A_2\succ A_5\succ A_4\succ A_3\succ A_1$ \\
        Shen et al.'s TOPSIS method in \cite{SMLXC2018} & $A_2\succ A_5\succ A_3\succ A_4\succ A_1$ \\
        Xu's IFWA operator method in \cite{Xu2007} & $A_2 \succ A_4\succ A_5\succ A_3\succ A_1$ \\
		Our TOPSIS based on $\varrho^{(100)}$ & $A_2 \succ A_5\succ A_4\succ A_3\succ A_1$ \\
        Our TOPSIS based on $\tilde{\varrho}^{(100)}$ & $A_2 \succ A_5\succ A_4\succ A_3\succ A_1$ \\
        Our TOPSIS based on $\varrho_{_{K_{0.2}, K_{0.4}}}^{(100)}$ & $A_2 \succ A_5\succ A_4\succ A_3\succ A_1$ \\
        Our TOPSIS based on $\varrho_{_{K_{0.5}, K_{0.4}}}^{(100)}$ & $A_2 \succ A_5\succ A_4\succ A_3\succ A_1$ \\
        Our TOPSIS based on $\varrho_{_{K_{0.6}, K_{0.4}}}^{(100)}$ & $A_2 \succ A_5\succ A_4\succ A_3\succ A_1$ \\
        \bottomrule
	\end{tabular}
      }
\end{table}

To illustrate the detailed influence of the parameters $\lambda$ and $\gamma_1$ on the
decision-making results in Example~\ref{Exm-WuX} by using metrics $\varrho^{(\lambda)}$,
$\tilde{\varrho}^{(\lambda)}$, and $\varrho_{_{K_{\gamma_1}, K_{\gamma_2}}}^{(\lambda)}$,
the relative closeness degrees $\mathscr{C}_i$ of each alternative $A_{i}$ obtained by
$\varrho^{(\lambda)}$, $\tilde{\varrho}^{(\lambda)}$, and $\varrho_{_{K_{\gamma_1},
K_{\gamma_2}}}^{(\lambda)}$ are shown in Figs.~\ref{Fig-Exm-11} (a), (b), and (c)--(d), respectively.
As can be seen from Fig.~\ref{Fig-Exm-11}, the ranking orders of the alternatives
using different values of parameters $\lambda$ and $\gamma_1$ remain the same and
are stabilized , when the parameter $\lambda$ is large enough, and thus the preferences
of decision makers do not affect the ranking results in this example.
This indicates that our method is effective and stable.

\begin{figure}[h]
\centering
\subfigure[Relative closeness degrees of $A_{1}$--$A_{5}$ obtained by $\varrho^{(\lambda)}$]
{\scalebox{0.295}{\includegraphics[]{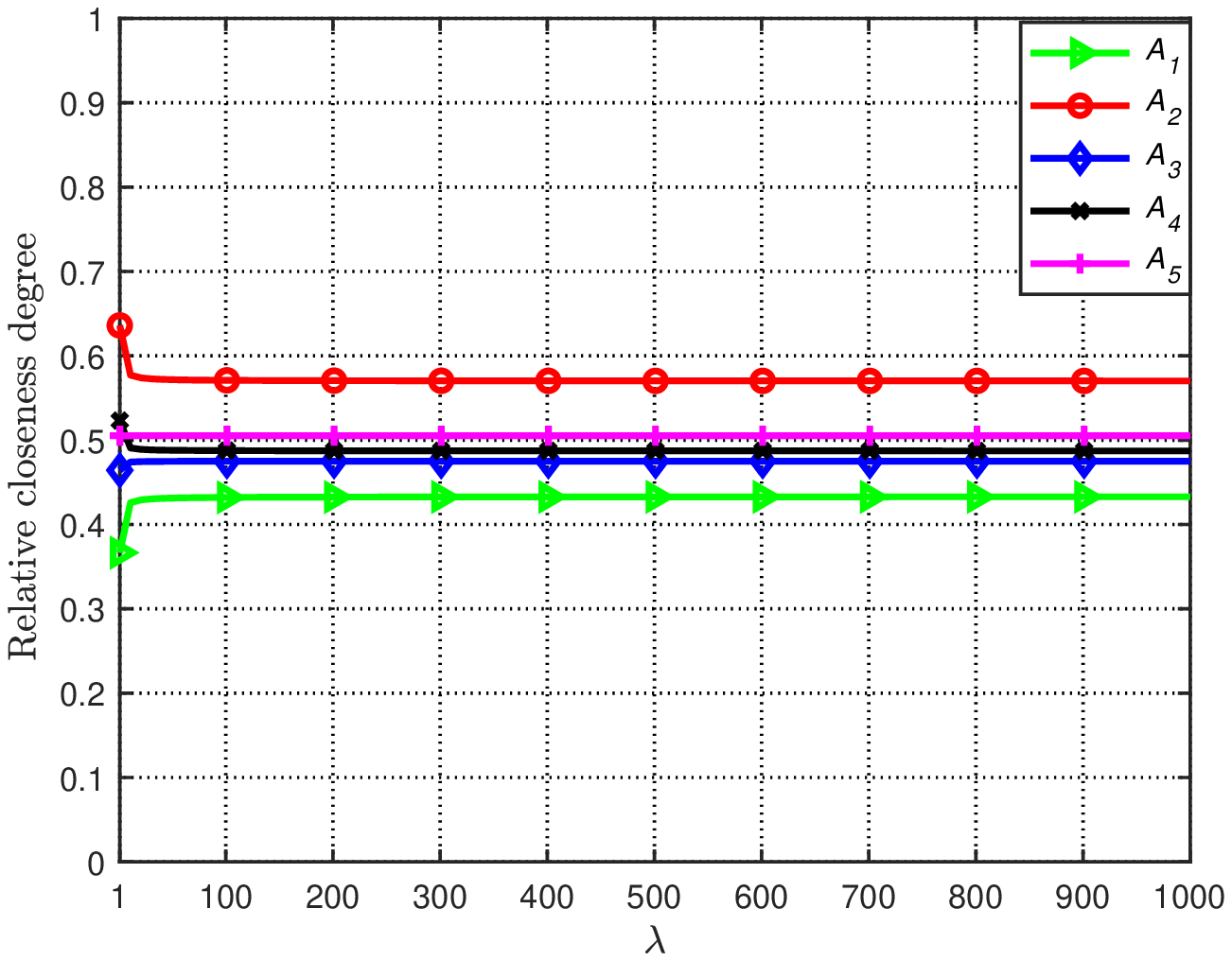}}}
\subfigure[Relative closeness degrees of $A_{1}$--$A_{5}$ obtained by $\tilde{\varrho}^{(\lambda)}$]
{\scalebox{0.295}{\includegraphics[]{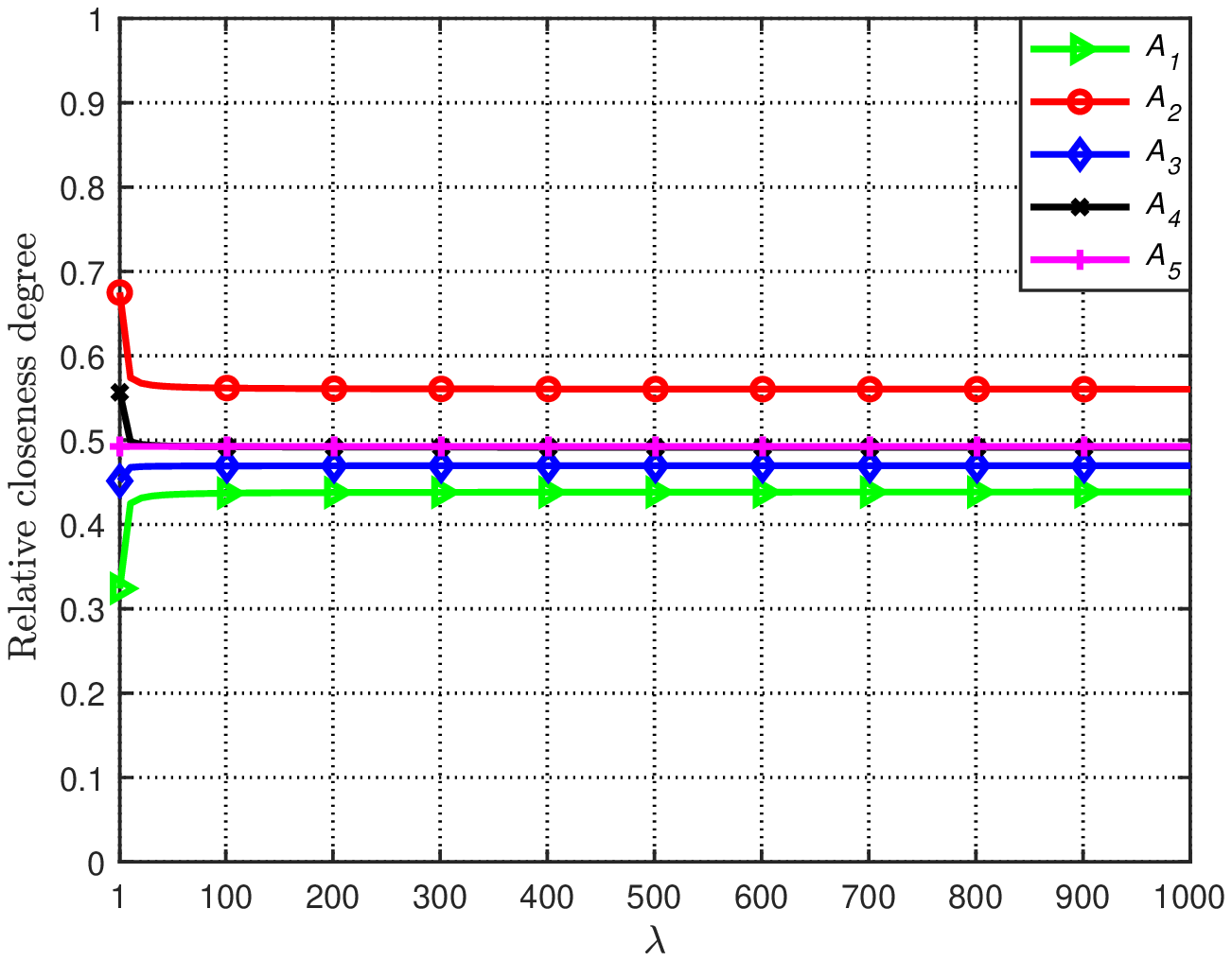}}}
\subfigure[Relative closeness degrees of $A_{1}$--$A_{5}$ obtained by $\varrho_{_{K_{\gamma_1},
K_{\gamma_2}}}^{(\lambda)}$ ($0\leq \gamma_1< 0.5$, $\gamma_2=1$)]
{\scalebox{0.295}{\includegraphics[]{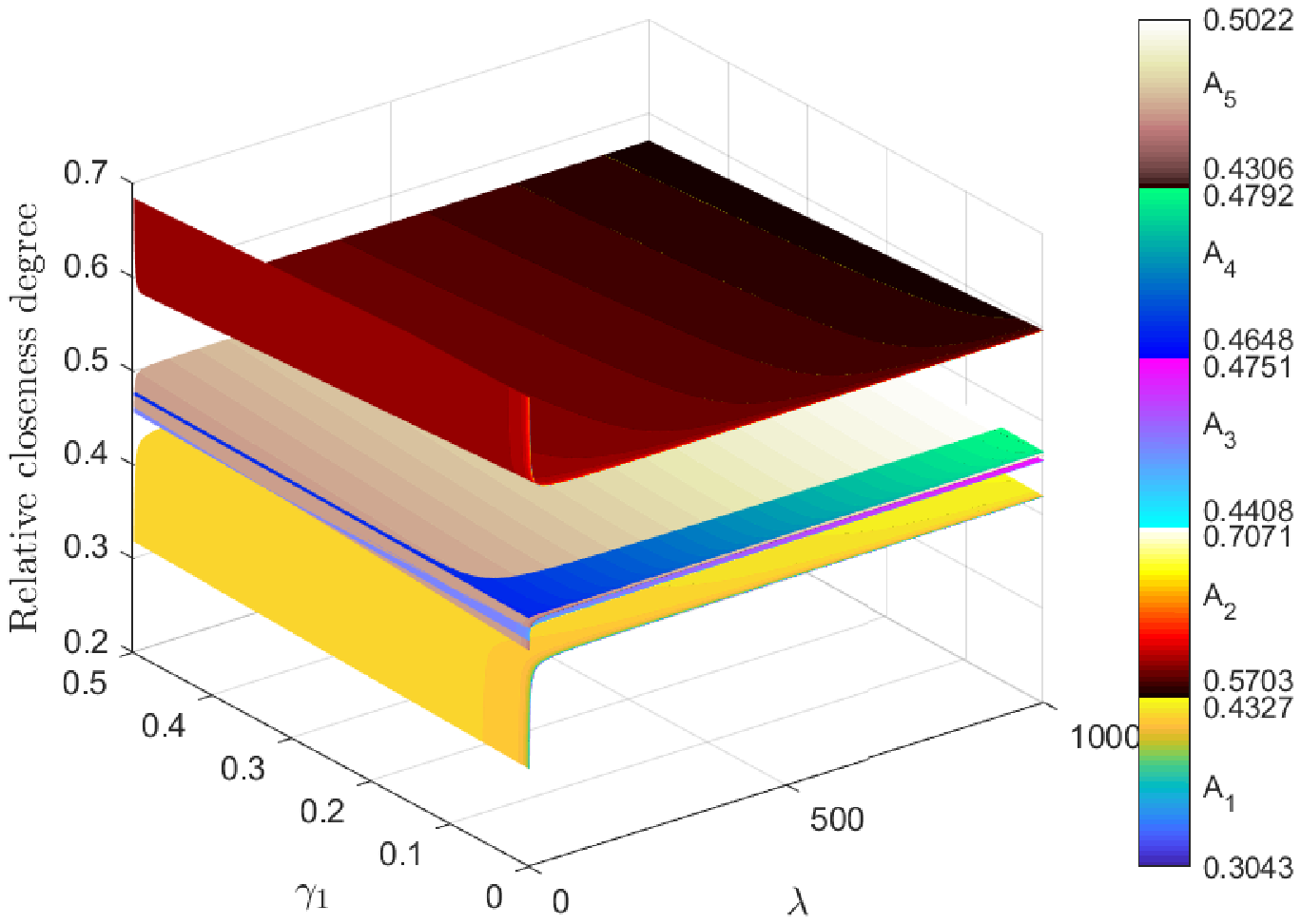}}}
\subfigure[Relative closeness degrees of $A_{1}$--$A_{5}$ obtained by $\varrho_{_{K_{\gamma_1},
K_{\gamma_2}}}^{(\lambda)}$ ($0.5< \gamma_1< 1$, $\gamma_2=1$)]
{\scalebox{0.295}{\includegraphics[]{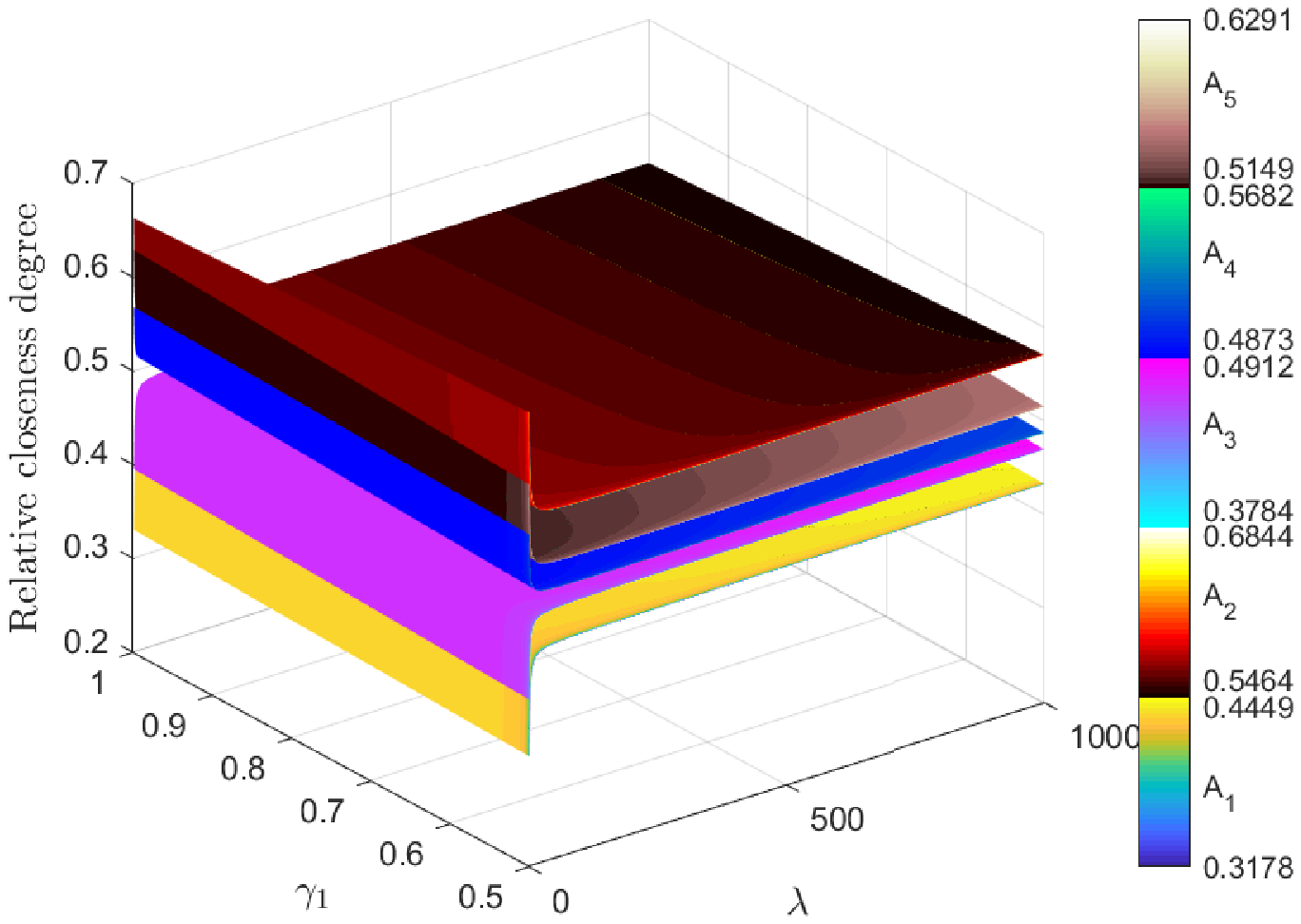}}}
\caption{Relative closeness degrees of $A_{1}$--$A_{5}$ in Example~\ref{Exm-WuX}}
\label{Fig-Exm-11}
\end{figure}

In summary, our proposed TOPSIS method has the following advantages:
\begin{enumerate}[(1)]
  \item It is monotonous under the linear order `$\leq_{_{\text{XY}}}$' or `$\leq_{_{\textrm{ZX}}}$' or
`$\leq_{_{A,B}}$'. This can overcome the limitation of non-monotonicity for some classical IF TOPSIS methods in~\cite{BBM2012,BGKA2009,CCL2016a,Li2014,MDJAR2019,RYU2020}. In addition, other IF TOPSIS methods can at
most guarantee the monotonicity under Atanassov's partial order $\subset$.
  \item Based on the admissible distances with linear orders, our method
is more in line with the essential features of the original TOPSIS introduced by Hwang and
Yoon~\cite{HY1981}.
  \item As can be seen from Example~\ref{Exm-WuX}, the ranking orders of the alternatives
using different MADM methods are slightly different. However, the best choice is the same.
Moreover, our preference order is stable when the parameter $\lambda$ is large enough.
This indicates that our method is effective and stable.
  \item Compared to the TOPSIS methods in \cite{SMLXC2018,CCL2016a}, our method requires less
  computation and fewer steps.
\end{enumerate}
}

\begin{example}
[{\textrm{\protect\cite[Example~5.2]{ZCK2019}}}]
\label{Exm-WuX-1}
Assume that there is a committee of a company, which decides to
choose a project manager from five alternatives $A_1$, $A_2$, $A_3$, $A_4$, and $A_5$
with four attributes $\mathscr{O}_1$, $\mathscr{O}_2$, $\mathscr{O}_3$, and $\mathscr{O}_4$,
where $\mathscr{O}_1$ is ``Self-Confidence", $\mathscr{O}_2$ is ``Personality",
$\mathscr{O}_3$ is ``Past Experience", $\mathscr{O}_4$ is the ``Proficiency in Project Management"
and $\mathscr{O}_1$, $\mathscr{O}_2$, $\mathscr{O}_3$, and $\mathscr{O}_4$ are all
benefit attributes, with weight vector $\omega=(0.1, 0.2, 0.3, 0.4)^{\top}$.

Assume the decision matrix $R=(r_{ij})_{5\times 4}$ given by the committee is as listed
in Table~\ref{Tab-3-Wu}.
\begin{table}[H]
\caption{The decision matrix $R$}
\label{Tab-3-Wu}\centering
\scalebox{0.85}{
	\begin{tabular}{cccccc}
		\toprule
		$$ & $\mathscr{O}_{1}$ &  $\mathscr{O}_{2}$ & $\mathscr{O}_{3}$ & $\mathscr{O}_{4}$ \\
		\midrule
		$A_{1}$ & $\langle0.4, 0.5\rangle$ & $\langle0.3, 0.6\rangle$ &
           $\langle0.4, 0.4\rangle$ & $\langle 0.5, 0.3\rangle$\\
		$A_{2}$ & $\langle0.4, 0.4\rangle$ & $\langle0.5, 0.4\rangle$ &
           $\langle0.3, 0.5\rangle$ & $\langle0.3, 0.4\rangle$ \\
        $A_{3}$ & $\langle0.4, 0.6\rangle$ & $\langle0.5, 0.5\rangle$ &
           $\langle0.4, 0.6\rangle$ & $\langle0.4, 0.6\rangle$ \\
        $A_{4}$ & $\langle0.3, 0.4\rangle$ & $\langle0.2, 0.6\rangle$ &
           $\langle0.1, 0.9\rangle$ & $\langle0.4, 0.4\rangle$ \\
        $A_{5}$ & $\langle0.5, 0.4\rangle$ & $\langle0.3, 0.6\rangle$ &
           $\langle0.3, 0.5\rangle$ & $\langle0.47, 0.5\rangle$ \\
        \bottomrule
	\end{tabular}
      }
\end{table}

Step~1: (Normalize the decision matrix) Since $\mathscr{O}_1$,
  $\mathscr{O}_2$, $\mathscr{O}_3$, and $\mathscr{O}_4$ are
all benefit attributes, we have $\overline{R}=(\bar{r}_{ij})_{5\times 4}=R$.

Step~2: (Determine the positive and negative ideal-points)
  The IF positive ideal-point is
  $$
  \mathbf{A}^{+}=(\langle 0.5, 0.4\rangle,
  \langle 0.5, 0.4\rangle, \langle 0.4, 0.4\rangle, \langle 0.5, 0.3\rangle)^{\top},
  $$
  and IF negative ideal-point is
  $$
  \mathbf{A}^{-}=(\langle 0.3, 0.6\rangle,
  \langle 0.2, 0.6\rangle, \langle 0.1, 0.9\rangle, \langle 0.3, 0.6\rangle)^{\top}.
  $$

  Step~3: (Compute the relative closeness degrees)
  Choose $\lambda=100$ and calculate the relative closeness degrees $\mathscr{C}_{i}$ of
  the alternatives $A_i$ ($i=1, 2, 3, 4, 5$) to the IF positive
  ideal-point $\mathbf{A}^+$ by Eqs.~\eqref{Eq-1}, \eqref{Eq-2}, and
  \eqref{Eq-3}:
  $\mathscr{C}_{1}=0.5565$, $\mathscr{C}_{2}=0.5295$,
  $\mathscr{C}_{3}=0.5058$, $\mathscr{C}_{4}=0.4555$,
  $\mathscr{C}_{5}=0.5171.$

  Step~4: (Rank the alternatives) Because
  $\mathscr{C}_1> \mathscr{C}_2> \mathscr{C}_5> \mathscr{C}_3
  > \mathscr{C}_4$, the ranking order of the alternatives $A_i$ ($i=1, 2, 3, 4, 5$) is:
   $A_1\succ A_2\succ A_5\succ A_3\succ A_4.$

Repeating Steps~1--2, by applying Eqs.~\eqref{Eq-1-a} and \eqref{Eq-2-b}, we obtain the following result:

Step~3: (Compute the relative closeness degrees)
  Calculate the relative closeness degrees $\mathscr{C}_{i}$ of
  the alternatives $A_i$ ($i=1, 2, 3, 4, 5$) to the IF positive
  ideal-point $\mathbf{A}^+$ by Eqs.~\eqref{Eq-1-a}, \eqref{Eq-2-b}, and
  \eqref{Eq-3}:
  $\mathscr{C}_{1}=0.5531$, $\mathscr{C}_{2}=0.5329$,
  $\mathscr{C}_{3}=0.5042$, $\mathscr{C}_{4}=0.4583$,
  $\mathscr{C}_{5}=0.5198.$

  Step~4: (Rank the alternative) Because
  $\mathscr{C}_1> \mathscr{C}_2> \mathscr{C}_5> \mathscr{C}_3
  > \mathscr{C}_4$, the ranking order of the alternatives $A_i$ ($i=1, 2, 3, 4, 5$) is:
   $A_1\succ A_2\succ A_5\succ A_3\succ A_4.$

{From Table~\ref{Tab-Exm-WuX-1}, which shows a comparison of the ranking orders of the alternatives
in Example~\ref{Exm-WuX-1} for different MADM methods, it can be observed that our results based on the metrics
$\varrho^{(100)}$, $\tilde{\varrho}^{(100)}$, $\varrho_{_{K_{0.5}, K_{0.4}}}^{(100)}$,
and $\varrho_{_{K_{0.6}, K_{0.4}}}^{(100)}$ are consistent
with the ranking orders obtained by the MADM methods
in~\cite{CCL2016a,WW2008,AAB2021,BG2016,ZZY2020,ZCK2019}.
\begin{table}[H]	
	\centering
	\caption{A comparison of the ranking orders of the alternatives
in Example~\ref{Exm-WuX-1} for different MADM methods}
	\label{Tab-Exm-WuX-1}
     \scalebox{0.85}{
	\begin{tabular}{cccccc}
		\toprule
		Methods & Ranking orders \\
		\midrule
		Chen et al.'s TOPSIS method in \cite{CCL2016a} & $A_1\succ A_2\succ A_5\succ A_3\succ A_4$ \\
		Wang and Wei's TOPSIS method in \cite{WW2008} & $A_1\succ A_2\succ A_5\succ A_3\succ A_4$ \\
        Altan~Koyuncu et al.'s TOPSIS method in \cite{AAB2021} & $A_1\succ A_2\succ A_5\succ A_3\succ A_4$ \\
        B{\"u}y{\"u}k{\"o}zkan and G{\"u}lery{\"u}z's TOPSIS method in \cite{BG2016} & $A_1\succ A_2\succ A_5\succ A_3\succ A_4$ \\
        Zhang et al.'s TOPSIS method in \cite{ZZY2020} & $A_1\succ A_2\succ A_5\succ A_3\succ A_4$ \\
        Zeng et al.'s VIKOR method in \cite{ZCK2019} & $A_1\succ A_2\succ A_5\succ A_3\succ A_4$ \\
		Our TOPSIS based on $\varrho^{(100)}$ & $A_1\succ A_2\succ A_5\succ A_3\succ A_4$ \\
        Our TOPSIS based on $\tilde{\varrho}^{(100)}$ & $A_1\succ A_2\succ A_5\succ A_3\succ A_4$ \\
        Our TOPSIS based on $\varrho_{_{K_{0.2}, K_{0.4}}}^{(100)}$ & $A_1\succ A_3\succ A_5\succ A_2\succ A_4$ \\
        Our TOPSIS based on $\varrho_{_{K_{0.5}, K_{0.4}}}^{(100)}$ & $A_1\succ A_2\succ A_5\succ A_3\succ A_4$ \\
        Our TOPSIS based on $\varrho_{_{K_{0.6}, K_{0.4}}}^{(100)}$ & $A_1\succ A_2\succ A_5\succ A_3\succ A_4$ \\
        \bottomrule
	\end{tabular}
      }
\end{table}}
\end{example}

\begin{figure}[h]
\centering
\subfigure[Relative closeness degrees of $A_{1}$--$A_{5}$ obtained by $\varrho^{(\lambda)}$]
{\scalebox{0.295}{\includegraphics[]{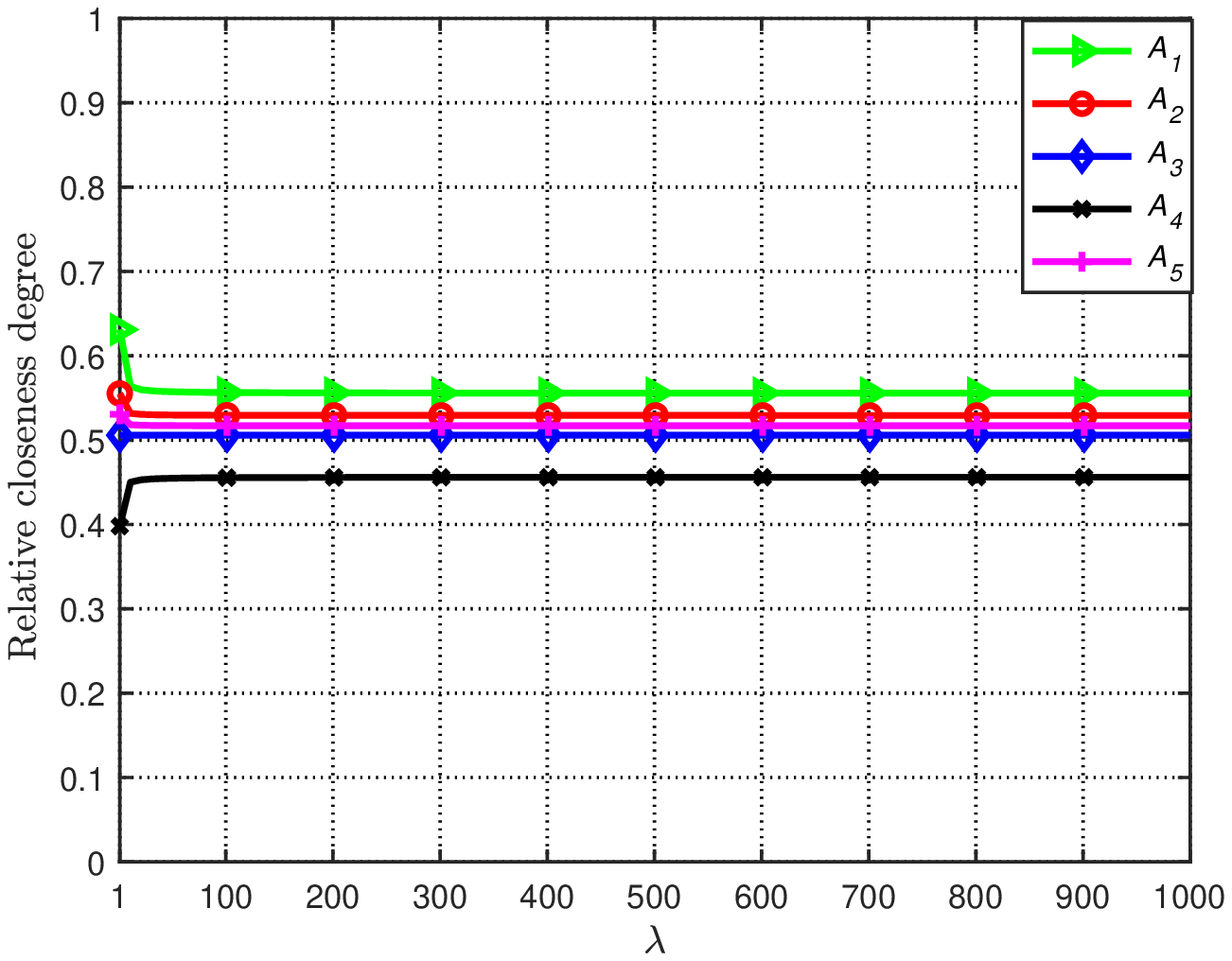}}}
\subfigure[Relative closeness degrees of $A_{1}$--$A_{5}$ obtained by $\tilde{\varrho}^{(\lambda)}$]
{\scalebox{0.295}{\includegraphics[]{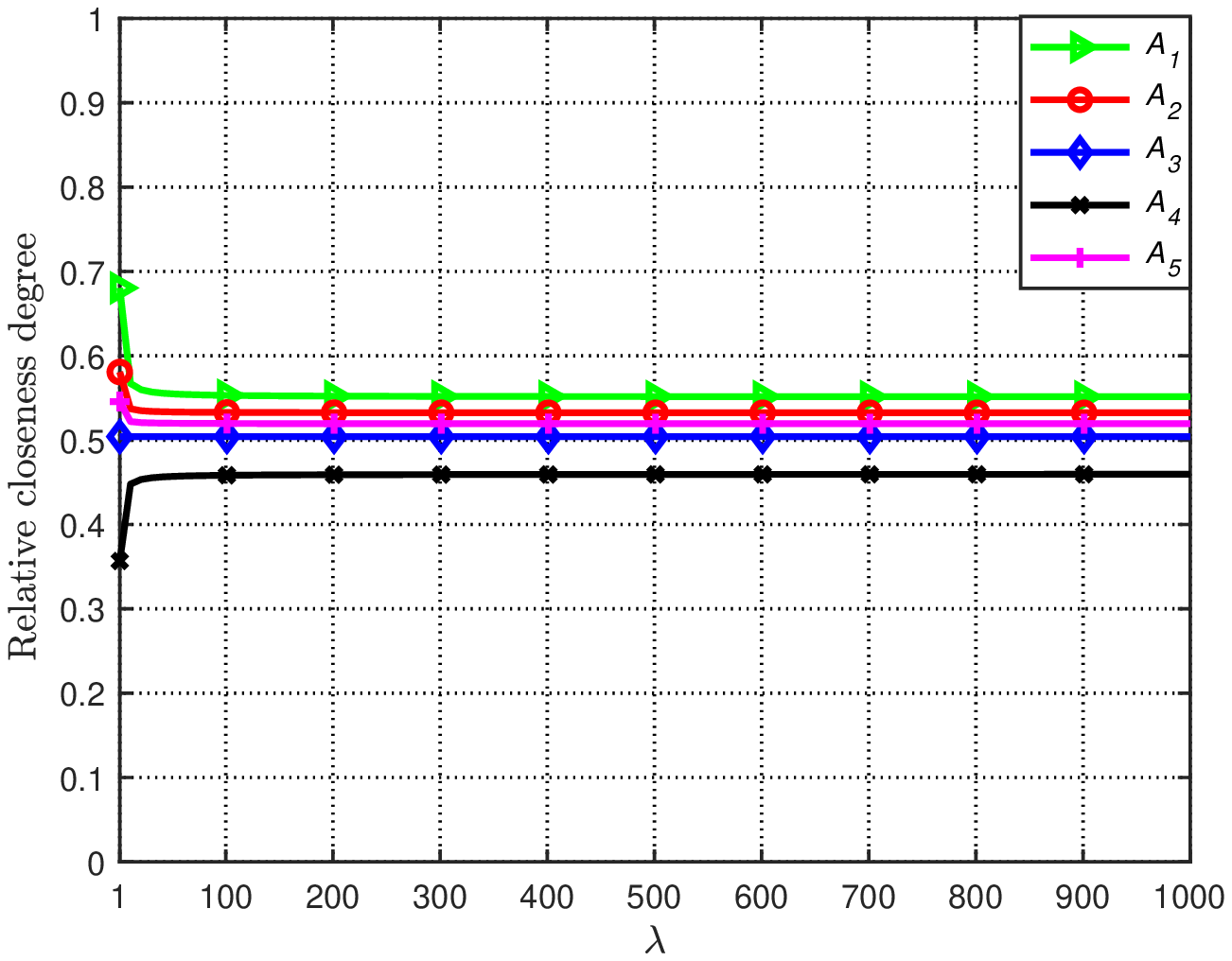}}}
\subfigure[Relative closeness degrees of $A_{1}$--$A_{5}$ obtained by $\varrho_{_{K_{\gamma_1},
K_{\gamma_2}}}^{(\lambda)}$ ($0\leq \gamma_1<0.5$, $\gamma_2=1$)]
{\scalebox{0.295}{\includegraphics[]{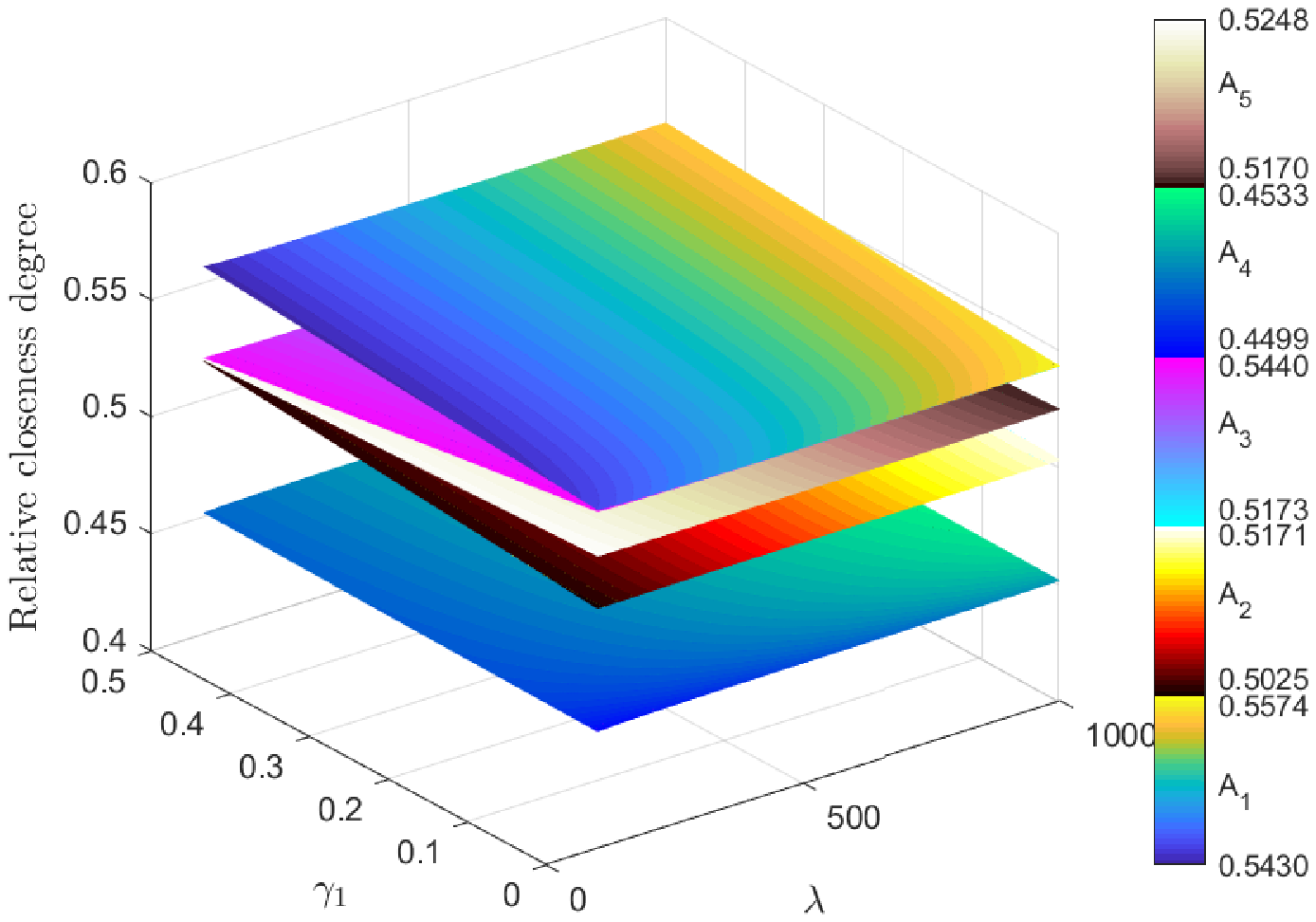}}}
\subfigure[Relative closeness degrees of $A_{1}$--$A_{5}$ obtained by $\varrho_{_{K_{\gamma_1},
K_{\gamma_2}}}^{(\lambda)}$ ($0.5< \gamma_1< 1$, $\gamma_2=1$)]
{\scalebox{0.295}{\includegraphics[]{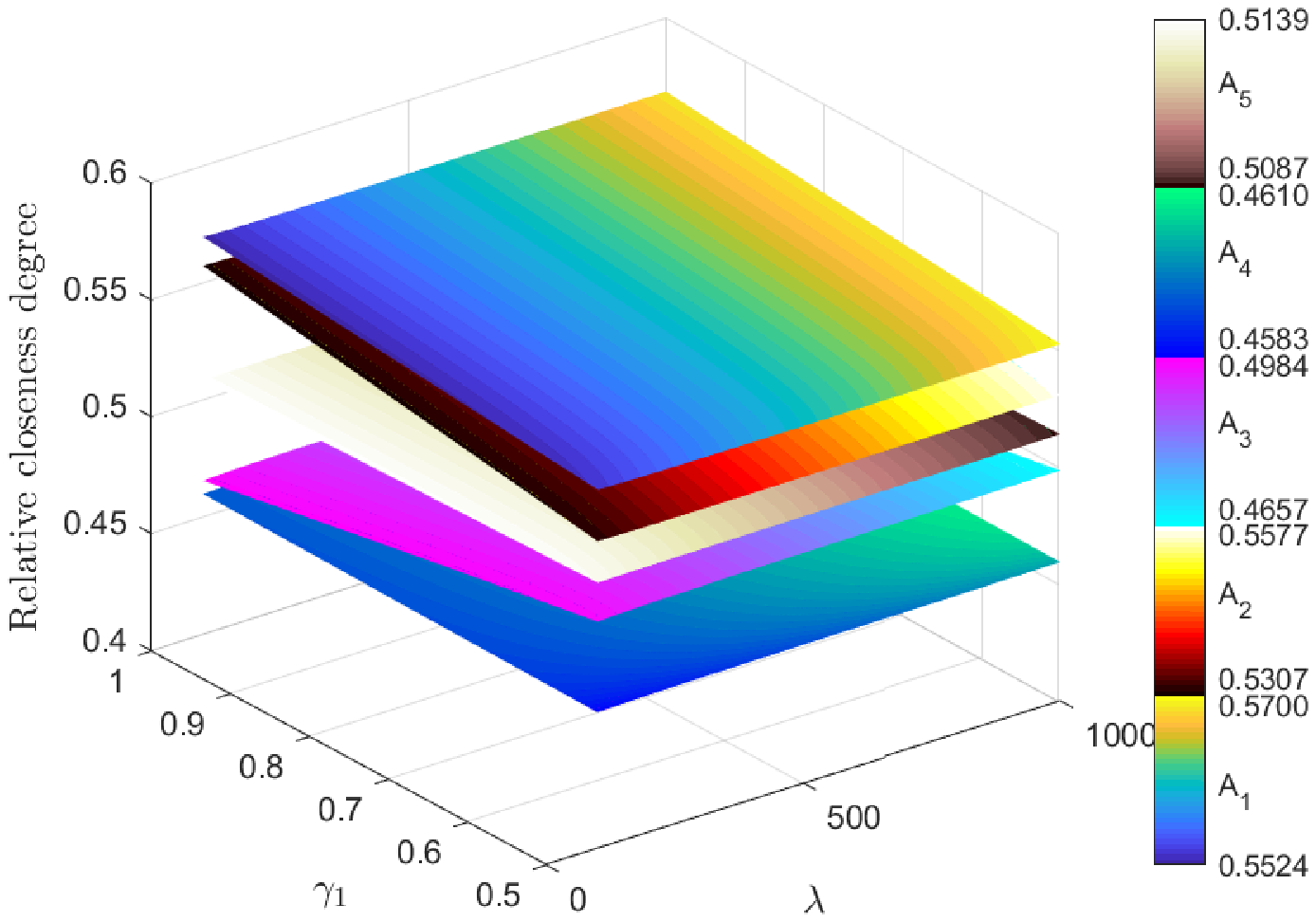}}}
\caption{Relative closeness degrees of $A_{1}$--$A_{5}$ in Example~\ref{Exm-WuX-1}}
\label{Fig-Exm-12}
\end{figure}

{To illustrate the detailed influence of the parameters $\lambda$ and $\gamma_1$ on the
decision-making results in Example~\ref{Exm-WuX-1} by using the metrics $\varrho^{(\lambda)}$,
$\tilde{\varrho}^{(\lambda)}$, and $\varrho_{_{K_{\gamma_1}, K_{\gamma_2}}}^{(\lambda)}$,
the relative closeness degrees $\mathscr{C}_i$ of each alternative $A_{i}$ obtained by
$\varrho^{(\lambda)}$, $\tilde{\varrho}^{(\lambda)}$, and $\varrho_{_{K_{\gamma_1},
K_{\gamma_2}}}^{(\lambda)}$ are shown in Fig.~\ref{Fig-Exm-12} (a), (b), and (c)--(d), respectively.
Compared to Example~\ref{Exm-WuX}, the preferences of decision makers greatly affect
the ranking results in this example: (1) For $0\leq \gamma_1<0.5$, the ranking order is:
$A_1\succ A_3\succ A_5\succ A_2\succ A_4$. (2) For $0.5\leq \gamma_1<1$, the ranking order is:
$A_1\succ A_2\succ A_5\succ A_3\succ A_4$.}

\section{Conclusions}\label{Sec-8}

{This paper is devoted to establishing a monotonous IF TOPSIS method with three typical linear orders,
`$\leq_{_{\text{XY}}}$' in~\cite{XY2006}, `$\leq_{_{\textrm{ZX}}}$' in~\cite{ZX2012}, and `$\leq_{_{A,B}}$' in \cite{DeBFIKM2016,DeBPBDaBMO2016}.} Noting that
the TOPSIS method is closely related to the order structure and the metric/similarity measure,
we first discuss some examples to show that some classical similarity measures in
\cite{Li2014,Sz2014,Xu2007a,XC2008}, including Euclidean similarity measure, Minkowski similarity
measure, and modified Euclidean similarity measure, do not satisfy the axiomatic definition of IF
similarity measures. Then, we prove the nonexistence of a continuous function that can
distinguish IFV by a real number and is increasing with Atanassov's order `$\subset$'. As
a direct corollary, we prove that there is no any continuous similarity measure that can
distinguish between each pair of IFVs. Moreover, we show some illustrative examples to demonstrate that
some classical IF TOPSIS methods in~\cite{BBM2012,BGKA2009,CCL2016a,Li2014,MDJAR2019,RYU2020}
are not monotonous with Atanassov's partial order `$\subset$' or the linear order
`$\leq_{_{\text{XY}}}$', which may yield counter-intuitive results. {To overcome
this limitation, by using three new parametric admissible distances with the linear order
`$\leq_{_{\text{XY}}}$' or `$\leq_{_{\textrm{ZX}}}$' or `$\leq_{_{A,B}}$',} we develop a novel IF TOPSIS method and prove
that it is monotonically increasing with these two linear orders. Finally, we show two practical
examples with comparative analysis to other MADM methods to illustrate the efficiency of our TOPSIS
method.

{Because the proposed TOPSIS method depends on the choice of the linear orders,
choosing an appropriate order for a given problem is very important for practical
decision-making. Meanwhile, because the construction method of admissible
distances with linear orders presented in this paper is relatively rough, which
fails to capture all properties of the the corresponding linear order, this may
cause inaccurate decision-making results in some cases. In the future, therefore,
we will further study the general construction of linear orders and admissible
distance/similarity measures for IFVs, which will be useful for building more effective
IF TOPSIS methods.}

\bibliographystyle{IEEEtran}
\bibliography{IEEEexample}


\end{document}